\documentclass[10pt, leqno]{amsart}
\usepackage{amscd}
\usepackage{amssymb}
\usepackage{amsmath} 
\usepackage[all]{xypic}
\usepackage{amsthm} 
\usepackage{stmaryrd}
\usepackage{mathabx}

\usepackage[usenames, dvipsnames]{color}

\theoremstyle{plain}
\newtheorem{Lem}{Lemma}[section]

\newtheorem{Cor}[Lem]{Corollary}
\newtheorem{Thm}[Lem]{Theorem}

{\theoremstyle{definition} 

\newtheorem{Ex}[Lem]{Example}
\newtheorem{Rk}[Lem]{Remark}
\newtheorem{Def}[Lem]{Definition}}

{\theoremstyle{definition} 

}

\newcommand{\zig}{\addtocounter{Lem}{1}\tag{\theLem}} 
\def\:{\colon} 
\DeclareMathOperator*{\colim}{colim}
\DeclareMathOperator*{\holim}{holim}
\DeclareMathOperator*{\limone}{lim^1}
\DeclareMathOperator*{\lims}{lim^\mathit{s}}
\DeclareMathOperator*{\hocolim}{hocolim}

\begin{document}
\title{A homotopy 
orbit spectrum for profinite groups}
\author[Daniel G. Davis, Vojislav Petrovi\'{c}]{Daniel G. Davis$\sp 1$, Vojislav Petrovi\'{c}}
\footnotetext[1]{The initial version of this paper \cite{orbitsone} was done while the first author was partially supported by a VIGRE NSF grant of the Purdue University Mathematics Department.}
\begin{abstract}
For a profinite group $G$, we define an $S[[G]]$-module to be a certain type of 
$G$-spectrum $X$ built from an inverse system $\{X_i\}_i$ of $G$-spectra, with each $X_i$ naturally a $G/N_i$-spectrum, where $N_i$ is an open normal subgroup  and $G \cong \lim_i G/N_i$. We define the 
homotopy orbit spectrum $X_{hG}$ and its homotopy orbit spectral sequence. We give results about when its $E_2$-term satisfies 
$E_2^{p,q} \cong \lim_i H_p(G/N_i, \pi_q(X_i))$. 
Our main result is that this occurs if $\{\pi_\ast(X_i)\}_i$ degreewise consists of compact Hausdorff abelian groups and continuous homomorphisms, with each $G/N_i$ acting continuously on $\pi_q(X_i)$ for all $q$. If $\pi_q(X_i)$ is additionally always profinite, 
then the $E_2$-term is the continuous homology of $G$ with coefficients in the graded profinite 
$\widehat{\mathbb{Z}}[[G]]$-module $\pi_\ast(X)$. Other results include theorems about Eilenberg-Mac Lane spectra and about when 
homotopy orbits preserve weak equivalences. 
\end{abstract}
\maketitle
\begin{section}{Introduction}
\par
Let $G$ be a finite group and let $X$ be a (left, naive) $G$-spectrum. 
Then the {\it homotopy orbit spectrum} $X_{hG}$ is defined to be 
$\hocolim_G X$, the homotopy colimit of the $G$-action on $X$ (see, 
for example, \cite[page 42]{book}). Furthermore, 
there is a homotopy orbit spectral sequence
\[H_p(G,\pi_q(X)) \Longrightarrow \pi_{p+q}(X_{hG}),\] where the $E_2$-term is 
the group homology of $G$, with coefficients in the graded $G$-module 
$\pi_\ast(X)$ (\cite[Section 5.1]{Mitchell}). 
In this paper, under certain hypotheses, we extend these constructions to the case 
where $G$ is a profinite group (at the end of this section, we give a discussion of related work).
\par
After making a few comments about notation, we summarize the contents of this paper. 
We follow the convention that all of our spectra are in $\mathrm{Spt}$, the category of Bousfield-Friedlander 
spectra of simplicial sets. We use $(-)_\mathtt{f}$ 
to denote functorial fibrant replacement in the category of 
spectra: for any spectrum $Z$, there is a natural map $Z \rightarrow Z_\mathtt{f}$ that is 
a weak equivalence, with $Z_\mathtt{f}$ fibrant. Also, ``holim'' always 
denotes the version of the homotopy limit of spectra that is constructed 
levelwise in the category of simplicial sets, as defined in 
\cite{Bousfield/Kan} and \cite[5.6]{Thomason}.
\par
In Section 2, for a finite group $G$, we use the fact that the homotopy 
colimit of a diagram of pointed simplicial sets is the diagonal of the 
simplicial replacement of the diagram, to obtain an alternative formulation of the 
homotopy orbits. Then, in Section 3, we use this formulation to define the 
homotopy orbit spectrum for a profinite group $G$ and a 
certain type of $G$-spectrum $X$, the category of which we now define. We point out that the concept of 
$S[[G]]$-module in the following definition was essentially first formulated by Mark Behrens. 
\begin{Def}\label{module}
Given a profinite group $G$, let $\{N_i\} := \{N_i\}_i$, indexed by $\{i\}$, be a cofinal subcollection of all the open normal subgroups 
of $G$, so that the diagram $\{G/N_i\}$ of finite discrete groups is an inverse system 
and there is the canonical 
isomorphism $G \cong \lim_i G/N_i$ of topological groups. We fix the collection $\{N_i\}$. We define 
\[S[[G]] :=\holim_i (S[G/N_i])_\mathtt{f},\] where for each $i$, $S[G/N_i] := S^0 \wedge (G/N_i)_+$. 
Here, $S^0$ is the sphere spectrum and more detail about $S[G/N_i]$ is in 
Definition \ref{groupring}. Let \[\{X_i\}_i\] be an inverse system of $G$-spectra and $G$-equivariant maps indexed over $\{i\}$, 
such that for each $i$, the $G$-action on $X_i$ factors through $G/N_i$ (so that $X_i$ is a 
$G/N_i$-spectrum) and $X_i$ is a fibrant spectrum. Then there is the 
$G$-spectrum \[X = \holim_i X_i,\] and the pair $(\{X_i\}_i, \holim_i X_i)$ is an {\em $S[[G]]$-module}. 
Such pairs are the objects of the {\em category of $S[[G]]$-modules}. Let $(\{Y_i\}_i, \holim_i Y_i)$ be 
an $S[[G]]$-module and let 
\[\tau \: \{X_i\}_i \to \{Y_i\}_i\] be a natural transformation of diagrams of $G$-spectra, where for 
each $i$, there is the component $\tau_i \: X_i \to Y_i$ of $\tau$. There is the induced $G$-equivariant map
\[\holim_i \tau_i \: \holim_i X_i \to \holim_i Y_i\] of spectra, and the pair $(\tau, \holim_i \tau_i)$ is a {\em morphism of $S[[G]]$-modules} \[(\{X_i\}_i, \holim_i X_i) \to (\{Y_i\}_i, \holim_i Y_i).\] 
Definition \ref{strong} defines composition of morphisms, completing the definition of the category of $S[[G]]$-modules. This category depends 
on the fixed family $\{N_i\}$, but we do not mention $\{N_i\}$ explicitly in the term ``category of $S[[G]]$-modules." 
\end{Def}

\begin{Ex}
In Definition \ref{module}, let $G$ act on each $G/N_i$ in the natural way and trivially on 
$S^0$, so that each $S[G/N_i]$ is a $G/N_i$-spectrum. Then so is each $(S[G/N_i])_\mathtt{f}$, by 
functoriality, and hence, the pair $(\{(S[G/N_i])_\mathtt{f}\}_i, S[[G]])$ is an $S[[G]]$-module. 
\end{Ex} 

\begin{Rk} 
Given a profinite group $G$, the expression ``let $X = \holim_i X_i$ be an $S[[G]]$-module" and natural variations of this expression mean 
that as in Definition \ref{module}, there is a fixed family $\{N_i\}_i$ 
and an $S[[G]]$-module $(\{X_i\}_i, \holim_i X_i)$, with 
$X = \holim_i X_i$.
\end{Rk} 

\par
After defining the homotopy orbit spectrum $X_{hG}$ for an $S[[G]]$-module $X = \holim_i X_i$, we 
show in Remarks \ref{reduction} and \ref{whenfinite} that this construction agrees with 
the classical definition when $G$ is finite. In Theorem \ref{generalized}, we prove that given a 
morphism $(\tau, \holim_i \tau_i)$ of $S[[G]]$-modules, as in Definition \ref{module}, if 
$\holim_i \tau_i$ is a weak equivalence of spectra, then the induced map
\[(\holim_i X_i)_{hG} \to (\holim_i Y_i)_{hG}\] is a weak equivalence. Here, when it is additionally true that the map 
$\lim_i \pi_\ast(\tau_i)$ is a bijection, Corollary \ref{3situations} describes three scenarios that 
imply that $\holim_i \tau_i$ is a weak equivalence. 
\par
We need the following notation. Let $G$ be any profinite group and let $\widehat{\mathbb{Z}}$ denote $\lim_{n \geq 1} \mathbb{Z}/n\mathbb{Z}$, the profinite 
completion of the integers and a commutative profinite ring. Then the complete group algebra $\widehat{\mathbb{Z}}[[G]]$ is defined by \[\widehat{\mathbb{Z}}[[G]] = \lim_{N \vartriangleleft_o G} \widehat{\mathbb{Z}}[G/N],\] where 
the inverse limit is over all the open normal subgroups of $G$ and each $\widehat{\mathbb{Z}}[G/N]$ is the group algebra. Let $B$ be a profinite left $\widehat{\mathbb{Z}}[[G]]$-module. Given a profinite right $\widehat{\mathbb{Z}}[[G]]$-module $A$, there is the complete tensor product $A \widehat{\otimes}_{\widehat{\mathbb{Z}}[[G]]} B$. By equipping $\widehat{\mathbb{Z}}$ with the trivial right $G$-action, we regard $\widehat{\mathbb{Z}}$ 
as a profinite right $\widehat{\mathbb{Z}}[[G]]$-module. When $p \geq 0$, let 
\[H^c_p(G,B) := \mathrm{Tor}_p^{\widehat{\mathbb{Z}}[[G]]}(\widehat{\mathbb{Z}},B)\] denote the $p$th continuous homology group of $G$, with coefficients in $B$. Here, the Tor group is defined in either of the usual two ways: for example, one can take the $p$th left derived functor of 
$(-)\widehat{\otimes}_{\widehat{\mathbb{Z}}[[G]]} B$ and apply it to $\widehat{\mathbb{Z}}$. For more information about the material in this paragraph, one helpful reference is \cite[Chapters 5, 6]{Ribes}. 

Now we state a result that consists of Theorems \ref{generalh-oSS} and \ref{corollary:whenfinite} and 
Corollary \ref{profinitecase}. 

\begin{Thm}\label{1stmain}
Let $G$ be any profinite group and let 
$X = \holim_i X_i$ be an $S[[G]]$-module. Then there is a 
homotopy orbit spectral sequence having the form
\[E_2^{p,q} 
\Longrightarrow \pi_{p+q}(X_{hG}),\] where $E_2^{p,q}$ is the $p$th homology of a certain Moore complex $\mathrm{(}$specified in Theorem \ref{generalh-oSS}$\mspace{1.5mu}\mathrm{)}$. If 
for every integer $t$, 
\begin{itemize}
\item
the inverse system $\{\pi_t(X_i)\}_i$ consists of compact Hausdorff abelian groups and continuous homomorphisms, and 
\item
for each $i$, the induced action of the discrete 
group $G/N_i$ on the compact Hausdorff space $\pi_t(X_i)$ is continuous,
\end{itemize} 
then there is the isomorphism
\begin{equation}\label{limitofhomologies}\zig
E_2^{p,q} \cong \lim_i H_p(G/N_i, \pi_q(X_i)).
\end{equation} 
Additionally, if for every $t$ and $i$, $\pi_t(X_i)$ is a profinite group, then 
\begin{equation}\label{obtainH^c}\zig
E_2^{p,q} \cong H^c_p(G,\pi_q(X)).
\end{equation}
\end{Thm} 

\begin{Rk}
The coefficients on the right-hand side of (\ref{obtainH^c}) deserve some explanation. 
In this situation, the two paragraphs after the proof of Theorem \ref{corollary:whenfinite} 
explain that for every integer $q$, the canonical $G$-equivariant map 
$\pi_q(X) \to \lim_i \pi_q(X_i)$ is an isomorphism of abelian groups and the target is a 
profinite $\widehat{\mathbb{Z}}[[G]]$-module. Thus, we identify the $G$-module $\pi_q(X)$ 
with $\lim_i \pi_q(X_i)$, and in this way, $\pi_q(X)$ is a 
profinite $\widehat{\mathbb{Z}}[[G]]$-module. 
\end{Rk}

\begin{Rk}
In the context of Theorem \ref{1stmain}, suppose that $M$ is a $G/N_i$-module for some 
$i$. If $M$ is finite, then when $M$ is equipped with the discrete topology, the 
$G/N_i$-action on $M$ is continuous. Thus, if $\holim_i X_i$ is an $S[[G]]$-module with 
$\pi_t(X_i)$ finite for every integer $t$ and all $i$, then we assign each $\pi_t(X_i)$ the discrete topology and the isomorphism in 
(\ref{obtainH^c}) holds.
\end{Rk}

To go further, we explain a condition that can be helpful to place on $G$. Recall that a topological space $Y$ is {\it 
countably based} if there is a countable family $\mathcal{B}$ of open sets, 
such that each open set of $Y$ is a union of members of $\mathcal{B}$. Then, if $G$ is a countably 
based profinite group, 
by \cite[Proposition 4.1.3]{Wilson}, $G$ has a chain 
\begin{equation}\label{chainofopens}\zig 
N_0 \geq N_1 \geq N_2 \geq \cdots
\end{equation}
of open normal subgroups, such that 
$G \cong \lim_{i\geq 0} G/N_i$. 

\begin{Ex}
Let $G$ be a compact $p$-adic analytic group. Then $G$ is a profinite group with an open 
subgroup $H$ of finite rank, by \cite[Corollary 8.34]{Dixon}. Hence, 
by \cite{Dixon}, $G$ has finite rank, and thus, is finitely 
generated, so that $G$ is a countably based profinite group. 
\end{Ex}

\begin{Def}\label{typeofmodule}
Let $G$ be a countably based profinite group and let $\{N_i\}_{i \geq 0}$ be a chain as in (\ref{chainofopens}), 
so that 
$G \cong \lim_{i \geq 0} G/N_i$. Also, suppose that 
\[X_0 \leftarrow X_1 \leftarrow 
\cdots \leftarrow X_i \leftarrow \cdots\] 
is a tower of $G$-spectra and $G$-equivariant maps such that with 
respect to 
the collection $\{N_i\}_{i \geq 0}$, the pair $(\{X_i\}_{i \geq 0}, \holim_i X_i)$ is 
an $S[[G]]$-module. Notice that each $X_i$ is a $G/N_i$-spectrum in the way prescribed by Definition \ref{module}. Then the pair 
$(\{X_i\}_{i \geq 0}, \holim_i X_i)$ is a {\em countably based 
$S[[G]]$-module}.
\end{Def} 

\begin{Rk}\label{rkforbased}
The expression ``let $X = \holim_i X_i$ be a countably based 
$S[[G]]$-module" and its natural variations mean that there is $G$ and a fixed chain $\{N_i\}_{i \geq 0}$, which are as in Definition \ref{typeofmodule}, and there is a countably based 
$S[[G]]$-module $(\{X_i\}_{i \geq 0}, \holim_i X_i)$, with $X = \holim_i X_i$. 
\end{Rk} 

When $K$ is a finite group and $Y$ is any spectrum equipped with the 
trivial $K$-action, the $K$-action on itself makes $Y[K]$ (see Definition \ref{groupring}) a $K$-spectrum, and there is an equivalence
\begin{equation}\label{classicalcontracting}\zig
(Y[K])_{hK} := \hocolim_K Y[K] \simeq Y
\end{equation} 
(for example, see \cite[page 127]{RognesStably}). 
Now let $G$ be a countably based profinite group. Given 
an arbitrary spectrum $Z$, in Definition \ref{Z[[G]]} we recall from chromatic homotopy theory the 
construction 
of $Z[[G]]$, an analogue of $Y[K]$, and we explain that $Z[[G]]$ and its underlying tower of $G$-spectra form 
a countably based $S[[G]]$-module. In Theorem \ref{Z[[G]]theorem}, we show that there 
is an equivalence
\[(Z[[G]])_{hG} \simeq Z,\] giving an analogue of (\ref{classicalcontracting}). For example, it 
follows that with respect to $\{N_i\}_{i \geq 0}$, as in Definition \ref{typeofmodule},
\[(S[[G]])_{hG} \simeq S^0.\] 

As suggested by Theorem \ref{1stmain}, we do not have, in general, a compact homological label for 
the $E_2$-term of the homotopy orbit spectral sequence. Thus, for the general situation, 
Theorem \ref{2bigconditions} isolates 
two key hypotheses that yield the isomorphism in (\ref{limitofhomologies}). 

Now suppose that $G$ is countably based and let 
$\holim_i X_i$ be a countably based $S[[G]]$-module. In this case, 
the situation with the $E_2$-term is better, and we prove in Theorem \ref{oss} that for any 
integer $q$, the groups $E_2^{\ast, q}$ fit into a long exact sequence in which the other terms are homology groups of 
complexes that are algebraically more meaningful than the complex for the general $E_2$-term. 
Theorem \ref{Hzero} shows that in this long exact sequence, 
the term to the left of $E_2^{0,q}$ is $\lim^1_i (\pi_{q+1}(X_i))_{\mspace{-1mu}{\scriptscriptstyle{G/N_i}}}\mspace{1mu}$, where 
each $(-)_{\scriptscriptstyle{G/N_i}}$ is the coinvariants functor. Now we fix $q$. 
There is a tower $\{\mathcal{C}(q)_i\}_i$ of non-negatively graded chain complexes, such that for each $i$, there are isomorphisms
\[H_\ast(\mathcal{C}(q)_i) \cong H_\ast(G/N_i, \pi_q(X_i)),\] where the 
left-hand side consists of the homology groups of $\mathcal{C}(q)_i$. Here, ``$\mathcal{C}(q)_i$" is 
abbreviated 
notation for a key construction that appears, for 
example, in (\ref{homology}) and (\ref{iso-of-complexes}). In Corollary \ref{mightneed}, we prove 
that if $\lim^1_i \pi_{q+1}(X_i) = 0$, then there is an isomorphism
\[E_2^{\ast,q} \xrightarrow{\,\cong\,} H_\ast(\lim_i \mathcal{C}(q)_i).\]

For each $l \geq 0$, let $\mathcal{C}(q)_{i,l}$ denote the $l$th group of chains of $\mathcal{C}(q)_i$. 
Theorem \ref{gottenwithML} and Remark \ref{finitelygenerated} imply that 
if 
\begin{itemize}
\item
the tower $\{\mathcal{C}(q)_{i,l}\}_i$ of abelian groups satisfies the Mittag-Leffler condition for each $l \geq 0$, 
\item
$\lim^1_i \pi_{q+1}(X_i) = 0$, and
\item 
for each $i$, $\pi_q(X_i)$ is a finitely generated $G/N_i$-module,
\end{itemize} 
then there is an isomorphism
\[E_2^{\ast,q} \xrightarrow{\cong} \lim_i H_\ast(G/N_i, \pi_q(X_i)).\] This last isomorphism is used to obtain Corollary \ref{E-Mf-g}, 
which is discussed briefly below. 


Now we give a definition that is an analogue for abelian groups of 
part of the concept of an $S[[G]]$-module. 

\begin{Def}\label{nicesystem} 
Let $G$ be any profinite group and, as in Definition \ref{module}, let $\{N_i\}$ be a fixed 
cofinal subcollection of 
all the open normal subgroups of $G$, indexed by a directed poset $\{i\}$. If $\{A_i\}$ is an inverse system 
in the category of $G$-modules indexed by $\{i\}$, such that for each 
$i$, the $G$-action on $A_i$ factors through $G/N_i$ (thus, $A_i$ is a 
$G/N_i$-module), then we call $\{A_i\}$ a {\em nice inverse system} of 
$G$-modules.
\end{Def}

When $A$ is an abelian group, we let $H(A)$ denote the Eilenberg-Mac Lane spectrum. 
Now let $\{A_i\}$ be a nice inverse system of $G$-modules. 
We show in Section \ref{EM} that 
the pair $(\{H(A_i)\}, \holim_i H(A_i))$ is an $S[[G]]$-module and by (\ref{e2term}), 
the $E_2$-term of the homotopy orbit spectral sequence, in general, can be simplified somewhat. 
Theorem \ref{generalE-Mresult} notes a condition that 
implies that 
\[\pi_\ast((\holim_i H(A_i))_{hG}) \cong \lim_i H_\ast(G/N_i,A_i),\] and when $G$ is 
countably based, Corollary \ref{E-Mf-g} gives 
hypotheses that yield this condition. This last result builds on a suggestion to the first author by Mark Behrens and his suggestion was the initial motivation for Section \ref{EM}. 
    
In \cite{PetrovicThesis}, the second author used the technique in the 
proof of Theorem \ref{oss} to study the $E_2$-term of a certain Adams-type 
spectral sequence for $\pi_\ast(L_{K(n)}Y)$, where $L_{K(n)}Y$ is the Bousfield localization with respect to the $n$th Morava $K$-theory $K(n)$ of an arbitrary spectrum $Y$. This technique again yields a long exact sequence and the second 
author studies its relationship to previous results in the literature about 
when this spectral sequence has $E_2$-term given by continuous cohomology.
  
\par
Others have investigated homotopy orbits for profinite 
groups. In \cite{fausk}, Halvard Fausk studies homotopy orbits in the setting of 
pro-orthogonal $G$-spectra, where $G$ is profinite. As 
pointed out in various places in this paper, Mark Behrens has worked on homotopy 
orbits for a profinite group. Also, Dan Isaksen has thought 
about homotopy orbits for profinite groups in the context of pro-spectra. 

After the appearance of \cite{orbitsone}, a manuscript that developed into the present paper, 
work on homotopy orbit 
spectra when $G$ is profinite was done by Gereon Quick in 
\cite[Section 2.5]{Quick}. Given a profinite $G$-spectrum 
$Y$, 
Quick defines a homotopy orbit spectrum $Y_{hG}$ 
and shows that there is a homotopy orbit spectral sequence 
\begin{equation}\label{QuickSpecSeq}\zig 
E^2_{s,t} = H_s^c(G; \pi_t(Y)) \Longrightarrow \pi_{s+t}(Y_{hG}),
\end{equation} where 
$E^2_{s,t}$ is given by continuous group homology. 
In 
\cite[Section 2.5: after Theorem 2.23]{Quick}, Quick discusses the relationship between these two works in the case when $G$ is countably based. The upshot is 
that in ``the rich world of profinite group actions and spectra $X$ with continuous action by such," sometimes 
$X$ is such that $X_{hG}$ and its homotopy orbit spectral sequence is given by \cite{Quick}, sometimes 
these are given by the present work, and 
there are instances when the constructions agree. We refer the reader to \cite{Quick} for more information about 
this situation. Also, if $G$ is any profinite group and $(\{X_i\}_i, \holim_i X_i)$ is an $S[[G]]$-module, 
each $X_i$ is a discrete $G$-spectrum, so that $\holim_i X_i$ need 
not be a profinite $G$-spectrum. Hence, there are $S[[G]]$-modules 
to which the construction of $(-)_{hG}$ in \cite{Quick} does not apply.

The paper \cite[Remark 1.8]{Szymik} contains 
an interesting remark about an application of a homotopy orbit spectrum with respect to $\mathbb{Z}_p$. The first author believes that the Lubin-Tate spectrum $E_n$ is an $S[[G_n]]$-module, where $G_n$ is the extended Morava stabilizer group, and he is working on this with Drew Heard. 

\vspace{.15in}
\par
\noindent
\textbf{Acknowledgements.}
The first author thanks Hal Sadofsky for sparking his interest in homotopy orbits for profinite 
groups, Mark Behrens for helpful and encouraging discussions 
about \cite{orbitsone} -- his positive influence occurs in multiple spots in the present paper, Drew Heard for 
discussions that partly motivated Theorem \ref{Z[[G]]theorem}, and 
Paul Goerss and Mark Hovey for their 
encouragement. We thank the referee of \cite{orbitsone}, the original version of this paper, for a helpful report, including a useful pointer about \cite[Theorem 5.3]{orbitsone} that helped the authors to improve upon it and eventually obtain Theorem \ref{corollary:whenfinite}. Also, we thank 
the referee of the revised versions for many helpful comments that improved the paper 
in a variety of ways. Finally, it was some unpublished explorations 
of Mike Hopkins and Hal Sadofsky about homotopy orbits for the action of the extended Morava stabilizer group on the Lubin-Tate spectrum that 
motivated the first author to begin the work on the project represented by this paper.
\end{section}
\begin{section}{Homotopy orbits when $G$ is finite}\label{two}
In this section, let $G$ be a finite group and let $X$ be a (left) $G$-spectrum. We use 
simplicial replacement to obtain a reformulation of the homotopy 
orbit spectrum $X_{hG} = \mathrm{hocolim}_G \,X$ that will be useful for defining the homotopy orbit 
spectrum when $G$ is a profinite group. Given a spectrum $Z$ and $k \geq 0$, we let $Z_k$ denote the
$k$th pointed simplicial set of $Z$ and, for $l \geq 0$, we use $Z_{k,l}$ to signify the set of 
$l$-simplices of $Z_k$. 
\par
Let $k \geq 0$ and let 
\[G \to X_k, \ \ \ \ast_G \mapsto X_k\] be the diagram defined by the action of $G$ on $X_k$, where, in the 
above display, $G$ is the one-object groupoid associated to the group and $\ast_G$ is the unique object in 
this groupoid. By \cite[Chapter XII, \S 5.2]{Bousfield/Kan}, 
\[(X_{hG})_k = \mathrm{hocolim}_G \, X_k = 
\mathrm{diag}(\, \textstyle{\coprod}_\ast \, (G \rightarrow X_k)),\] where 
$\mathrm{diag}(-)$ is the functor that takes the diagonal of a bisimplicial 
set, and $\coprod_\ast (G \rightarrow X_k)$ is the simplicial replacement of the diagram 
$G \rightarrow X_k$. In each dimension $l$, the simplicial
replacement is the 
pointed simplicial set
\[\textstyle{\coprod}_l \, 
(G \rightarrow X_k) := \textstyle{\bigvee}_{G^l}X_k,\] the wedge of copies of $X_k$ indexed 
by $G^l$, the $l$-fold product of copies of $G$. When $l = 0$, $G^l = \{e\}$, the trivial group, so that 
$\bigvee_{G^0} X_k = X_k$. 

Now we define the face and degeneracy maps of the simplicial replacement. Let $l \geq 0$. 
We identify the indexing 
set $G^l$ above with the set of all $l$-fold compositions 
\[
\ast_G \xleftarrow{g_1} \ast_G \xleftarrow{g_2} \cdots 
\xleftarrow{g_l} \ast_G,
\]
where $g_1, g_2, ..., g_l$, as morphisms in the groupoid $G$, are elements of 
the group $G$. Under this identification, $G^0 = \{e\}$ is regarded as the set $\{\ast_G\}$, where $\ast_G$ is a 
``$0$-fold composition." 

Suppose $l \geq 1$ and let $0 \leq j \leq l$. The $j$th face map 
$d_j \: \bigvee_{G^l} X_k \rightarrow \bigvee_{G^{l-1}} X_k$ is obtained from the universal property of the 
coproduct $\bigvee_{G^l} X_k$. Thus, to define $d_j$, for each $l$-fold composition
\[c_l := (\ast_G \overset{g_1}{\longleftarrow} \cdots \overset{g_l}{\longleftarrow} \ast_G),\] it suffices to define a map 
$d_j(c_l) \: X_k \to \bigvee_{G^{l-1}} X_k$, whose source is the copy of $X_k$ in $\bigvee_{G^l} X_k$ indexed by $c_l$. 
For the first two cases -- when $j = 0$ and $0 < j < l$, $d_j(c_l)$ is the inclusion into the coproduct 
of the copy of $X_k$ indexed by the $(l-1)$-fold composition 
\[\ast_G \overset{g_2}{\longleftarrow} \cdots 
\overset{g_l}{\longleftarrow} \ast_G\] and 
\[\ast_G \xleftarrow{g_1} \cdots \xleftarrow{g_{j-1}} \ast_G
\xleftarrow{g_j g_{j+1}} \ast_G \xleftarrow{g_{j+2}} \cdots
\xleftarrow{g_l} \ast_G,\] respectively. In the last case -- when $j = l$, $d_j(c_l)$ is the composition 
\[X_k \to X_k \to \textstyle{\bigvee}_{G^{l-1}} X_k,\] where the left-hand map $X_k \to X_k$ is the map in the diagram 
$G \to X_k$ determined by the morphism $\ast_G \xleftarrow{g_l} \ast_G$ (that is, $G \to X_k$ is a functor and $X_k \to X_k$ is its value on $g_l$) and the right-hand map is the inclusion of the copy of $X_k$ 
indexed by the $(l-1)$-fold composition $\ast_G \xleftarrow{g_1} \cdots 
\xleftarrow{g_{l-1}} \ast_G$. 

Let $l \geq 0$ and again, let $0 \leq j \leq l$. As above, the 
$j$th degeneracy map $s_j \: \bigvee_{G^l} X_k \rightarrow 
\bigvee_{G^{l+1}} X_k$ comes from the universal property of the coproduct $\bigvee_{G^l} X_k$, and so to define $s_j$, we only need to define maps $s_j(c_l) \: X_k \to \bigvee_{G^{l+1}} X_k$. Here, 
$c_l$ is an $l$-fold composition as before, $c_0$ denotes the unique $0$-fold composition $\ast_G$, and 
$s_j(c_l)$ has source equal to the copy of $X_k$ in the source of $s_j$ indexed by $c_l$. Then $s_j(c_l)$ is the 
inclusion into $\bigvee_{G^{l+1}} X_k$ of the copy of $X_k$ indexed by the $(l+1)$-fold composition
\[\ast_G \xleftarrow{g_1} \cdots \xleftarrow{g_{j}} \ast_G \xleftarrow{e} \ast_G \xleftarrow{g_{j+1}} \cdots 
\xleftarrow{g_l} \ast_G,\] where $e$ denotes the identity element of $G$. This completes our recollection of the 
definition of the simplicial replacement $\coprod_\ast (G \to X_k)$. 

\begin{Def}\label{groupring}
Let $K$ be a pointed simplicial set and let $L$ be a set. Then 
\[K[L]=K \wedge L_+,\] where $L_+$ is the constant simplicial set on $L$, 
together with a disjoint basepoint. Similarly, if $Z$ is a spectrum, then 
$Z[L]$ is the spectrum with each $(Z[L])_k$ equal to $Z_k \wedge L_+.$ We can also write 
$Z[L]$ as $Z \wedge L_+$ and we note that there is a natural isomorphism
\[Z[L] \cong \textstyle{\bigvee_L Z},\] where the right-hand side is a coproduct of copies of $Z$, indexed by the elements of $L$. 
\end{Def}
In dimension $l$, we have
\[\textstyle{\coprod}_l (G \rightarrow X_k) = 
\textstyle{\bigvee}_{G^l} X_k \cong X_k \wedge (G^l)_+ =
X_k[G^l].\] Since the middle isomorphism above is natural, there is a 
simplicial pointed simplicial set $X_k[G^\bullet]$ whose $l$-simplices are the 
pointed simplicial set $X_k[G^l]$ and which satisfies the isomorphism
\[X_k[G^\bullet] \cong \textstyle{\coprod_\ast (G \to X_k)}\] of simplicial pointed simplicial sets. It 
follows that 
\[(X_{hG})_k \cong \mathrm{diag}(X_k[G^\bullet]).\] 
\par
We introduce some notation that organizes the above, allowing us to summarize 
it in Theorem \ref{model} below. Though this result is labeled ``Theorem," it is just a repackaging of a standard definition in a way that is helpful to us in the next section (see Definition \ref{def}).
\begin{Def}\label{xgbullet}
Let $X$ be a $G$-spectrum. Then $X[G^\bullet]$ is the simplicial 
spectrum that is defined above, 
with $(X[G^\bullet])_k = X_k[G^\bullet]$ and $l$-simplices equal to 
the spectrum $X[G^l]$. Thus, $X[G^\bullet]$ is the simplicial 
spectrum
\[\xymatrix@C.25in{
X \cong X[\ast] & X[G] \ar@<.5ex>[l] \ar@<-.5ex>[l]
& X[G^2] \ar@<.7ex>[l] \ar[l] \ar@<-.7ex>[l]
& \cdots . \ar@<.8ex>[l] \ar@<.25ex>[l] \ar@<-.25ex>[l] \ar@<-.8ex>[l]}\] 
We let $\pi_\ast(X[G^\bullet])$ denote the simplicial (graded) abelian group associated to 
$X[G^\bullet]$.
\end{Def}
\begin{Def}
Let $d(Z_\bullet)$ denote the spectrum that is the 
diagonal of the simplicial spectrum $Z_\bullet$ (see \cite[page 100]{Jardine}). 
Thus, for all $k \geq 0,$ 
\[(d(Z_\bullet))_k = \mathrm{diag}((Z_\bullet)_k)\] and, for all $l \geq 0,$ 
\[(d(Z_\bullet))_{k,l} = (Z_l)_{k,l}.\] 
\end{Def}
\begin{Thm}\label{model}
If $G$ is a finite group and $X$ is a $G$-spectrum, then
\[X_{hG} \cong d(X[G^\bullet]).\] 
\end{Thm}

\begin{Rk}\label{bar}
In Theorem 2.4, by \cite[page 100]{Jardine}, $d(X[G^\bullet])$ can be viewed as 
the geometric realization (that is, as a type of coend) of the simplicial spectrum 
\[X[G^\bullet] = X \wedge (G^\bullet)_+,\] which is a bar construction for the action of $G$ on $X$. The notation ``$X[G^\bullet]$" reflects the algebra that arises when working with the Moore complex associated to 
the simplicial (graded) abelian group $\pi_\ast(X[G^\bullet])$. This is seen, for example, in the algebraic manipulations in the proof of Theorem \ref{corollary:whenfinite} and when $X$ is a homotopy ring spectrum, where for each $l \geq 0$, $X[G^l]$ is a ``homotopy group ring spectrum," with $\pi_0(X[G^l]) \cong \pi_0(X)[G^l]$. In this isomorphism, the right-hand side is a group ring, since $\pi_0(X)$ is a ring. 
\end{Rk}

\end{section}
\begin{section}{The homotopy orbit spectrum $X_{hG}$ when $G$ is profinite}\label{basicdefinition} 
\par
In this section, we use the model for the homotopy orbit spectrum that is given in 
Theorem \ref{model} to help us define homotopy orbits for $G$ a 
profinite group. We begin with a comment about Definition \ref{module}.
\begin{Rk}\label{remark}
Let $X = \holim_i X_i$ be an $S[[G]]$-module. The choice of the term 
``$S[[G]]$-module" is motivated by the fact that, for each $i$, the action 
of $G/N_i$ on $X_i$ yields a function $G/N_i \rightarrow 
\mathrm{Hom}_\mathrm{Spt}(X_i,X_i)$ in 
\begin{align*}
\mathrm{Hom}_\mathbf{Sets}(G/N_i,\mathrm{Hom}_\mathrm{Spt}(X_i,X_i)) 
& \cong \textstyle{\prod}_{G/N_i} \mathrm{Hom}_\mathrm{Spt}(X_i,X_i) \\ & \cong 
\mathrm{Hom}_\mathrm{Spt}(X_i[G/N_i],X_i)\end{align*} that corresponds to a map 
$X_i \wedge (G/N_i)_+ 
\rightarrow X_i$. The analogy with actual modules over the 
spectrum $S[[G]]$ could be pursued further, as Mark Behrens has done in unpublished work; however, for simplicity, we do not 
do this.
\end{Rk} 
\par
Using the first author's preliminary work on homotopy orbits as a reference point, the following definition was basically given by Behrens, and then later independently formulated by the first author.
\begin{Def}\label{def}
Let $G$ be any profinite group and let $(\{X_i\}_i, \holim_i X_i)$ be an $S[[G]]$-module, with 
$X = \holim_i X_i$. For each $l \geq 0$, the diagrams 
$\{X_i\} := \{X_i\}_i$ and $\{G/N_i\}$ induce the diagram \[\{X_i[(G/N_i)^l]\}\] 
in the following way: if 
$i' \leq i$, then (a) there are maps $\xi_{i,i'} \: X_i \to X_{i'}$ from $\{X_i\}$ and 
$\gamma_{i,i'} \: G/N_i \to G/N_{i'}$ from $\{G/N_i\}$; (b) $\gamma_{i,i'}$ induces a 
map $\gamma^l_{i,i'} \: (G/N_i)^l \to (G/N_{i'})^l$ in the natural way; (c) $\gamma^l_{i,i'}$ induces 
a map \[(\gamma^l_{i,i'})_+ \: ((G/N_i)^l)_+ \to ((G/N_{i'})^l)_+\] of pointed simplicial sets, 
where, for example, $((G/N_i)^l)_+$ is the constant simplicial set on $(G/N_i)^l$, together with a 
disjoint basepoint; and (d) thus, 
there is the map  
\[X_i[(G/N_i)^l] = X_i \wedge ((G/N_i)^l)_+ \xrightarrow{\,\xi_{i,i'} \wedge (\gamma^l_{i,i'})_+\,} X_{i'} \wedge ((G/N_{i'})^l)_+ = X_{i'}[(G/N_{i'})^l],\] which is the map in the diagram 
$\{X_i[(G/N_i)^l]\}$ associated to the relation $i' \leq i$. Given $l \geq 0$ and the diagram 
$\{X_i[(G/N_i)^l]\}$, one can 
form $\holim_i (X_i[(G/N_i)^l])_\mathtt{f}$, which gives the simplicial spectrum 
$\holim_i(X_i[(G/N_i)^\bullet])_\mathtt{f}$. Then 
we define $X_{hG}$, the {\em homotopy orbit 
spectrum} of the $S[[G]]$-module $X$ with respect to the $G$-action, to be the spectrum 
\[X_{hG}=d(\holim_i(X_i[(G/N_i)^\bullet])_\mathtt{f}).\] 
The functor $(-)_\mathtt{f}$ appears here so that the homotopy limit 
is well-behaved. 
\end{Def}

Henceforth, when $G$ is finite and $X$ is a $G$-spectrum, we use the notation \[X_{h'G} := \hocolim_G X\] to 
denote the classical homotopy orbit spectrum. In this case, in Remarks \ref{reduction} and \ref{whenfinite} 
below, we show that Definition \ref{def} reduces to the classical definition and, for any $G$-spectrum, the 
former definition recovers the latter one. 

\par
Definition \ref{def} comes from imitating the model in 
Theorem \ref{model} and from the demands of the homotopy orbit spectral sequence, especially its $E_2$-term, when 
it can be related to group homology. For example, the occurrences of ``$\holim_i$" 
in Definition \ref{def} correspond to the instrumental instances of ``$\lim_i$" 
that appear in the proof of Theorem \ref{corollary:whenfinite}. Also, since we use the spectral sequence in Theorem 
\ref{dss} below to build the homotopy orbit spectral sequence, we want $d(-)$ to be on the outside in the construction 
of $X_{hG}$, instead of on the inside, as in 
\begin{equation}\zig\label{justlimit}
\holim_i(d(X_i[(G/N_i)^\bullet]))_\mathtt{f}.
\end{equation} For each $i$, Theorem \ref{model} gives the isomorphism
\[(X_i)_{h'G/N_i} \cong d(X_i[(G/N_i)^\bullet]),\] so that there is a weak equivalence between the homotopy limit in (\ref{justlimit}) and $\holim_i ((X_i)_{h'G/N_i})_\mathtt{f}$. While the construction $\holim_i ((X_i)_{h'G/N_i})_\mathtt{f}$ is interesting, the natural 
spectral sequence associated to it is the homotopy spectral sequence
\[\lims_i \pi_t((X_i)_{h'G/N_i}) \Longrightarrow \pi_{t-s}(\holim_i ((X_i)_{h'G/N_i})_\mathtt{f}),\] which cannot yield a homotopy orbit spectral sequence that has $E_2$-term given by continuous group homology. 

The following comment was suggested by the referee of the revised versions of this paper. 

\begin{Rk}\label{bar2}
In the context of Definition \ref{def}, by Remark \ref{bar}, each $(X_i[(G/N_i)^\bullet])_\mathtt{f}$ is a bar 
construction for the action of $G/N_i$ on $X_i$, and thus, we can define 
\[B^\mathrm{cts}_\bullet(G,X) := \holim_i(X_i[(G/N_i)^\bullet])_\mathtt{f},\] so that we have the simple 
expression $X_{hG} = d(B^\mathrm{cts}_\bullet(G,X))$ and we can think of $B^\mathrm{cts}_\bullet(G,X)$ 
as being the ``continuous bar construction" for the action of $G$ on the $S[[G]]$-module $X$. 
\end{Rk}  

\par
We establish some notation and recall a useful result.

\begin{Def}\label{Moore}
If $X_\bullet$ is a simplicial spectrum, then for each integer $q$, we let $\pi_q(X_\ast)$ denote the Moore 
complex of the simplicial abelian group $\pi_q(X_\bullet)$. By ``Moore complex," we mean the chain complex 
\[\cdots \xrightarrow{\partial_2} \pi_q(X_2) \xrightarrow{\partial_1} \pi_q(X_1) \xrightarrow{\partial_0} \pi_q(X_0),\] with boundary homomorphism 
\[\partial_p = \textstyle{\sum}_{j=0}^{p+1} (-1)^j d_j \: \pi_q(X_{p+1}) \to \pi_q(X_{p}),\] for each non-negative integer $p$. Here, 
$d_0, d_1, ..., d_{p+1} \: \pi_q(X_{p+1}) \to \pi_q(X_{p})$ refer to face maps of $\pi_q(X_\bullet)$. Also, $H_p(\pi_q(X_\ast))$ denotes the $p$th homology of the Moore complex. In general, when our notation for a simplicial abelian group contains ``$\scriptstyle{\bullet}$," then we use the same notation for the Moore complex, but with ``$\scriptstyle{\bullet}$" changed to ``$\scriptstyle{\ast}$," as done above with $\pi_q(X_\ast)$. For example, in the context of Definition \ref{def}, for each $i$, $\pi_q((X_i[(G/N_i)^\ast])_\mathtt{f})$ is the Moore complex for the 
simplicial abelian group $\pi_q((X_i[(G/N_i)^\bullet])_\mathtt{f})$. 
\end{Def}

\begin{Thm}[{\cite[Corollary 4.22]{Jardine}}]\label{dss}
If $X_\bullet$ is a simplicial spectrum, then there is a spectral sequence
\[E_2^{p,q} = H_p(\pi_q(X_\ast)) \Longrightarrow \pi_{p+q}(d(X_\bullet)).\] 
\end{Thm}

As in Remark \ref{bar}, the abutment of the above spectral sequence can be viewed as the homotopy groups of the geometric realization of $X_\bullet$, which is a common way to understand this spectral sequence.

\begin{Lem}\label{we}
Let $X_\bullet \rightarrow Y_\bullet$ be a map between simplicial spectra, such that for each $n \geq 0$, the map $X_n \rightarrow Y_n$ is a weak equivalence between the 
$n$-simplices. Then the induced 
map $d(X_\bullet) \rightarrow d(Y_\bullet)$ is a weak equivalence of spectra.
\end{Lem}
\begin{proof}
There is a spectral sequence \[H_p(\pi_q(X_\ast)) \Longrightarrow \pi_{p+q}(d(X_\bullet)),\] and 
a map to the spectral sequence 
\[H_p(\pi_q(Y_\ast)) \Longrightarrow \pi_{p+q}(d(Y_\bullet)).\] Since $\pi_q(X_n) \cong \pi_q(Y_n)$, for each $n \geq 0$, and $\pi_q(X_\ast)$ and $\pi_q(Y_\ast)$ are chain complexes, there is an isomorphism $H_p(\pi_q(X_\ast)) \xrightarrow{\,\cong\,} H_p(\pi_q(Y_\ast))$ of $E_2$-terms. Therefore, the abutments of the above two spectral sequences are isomorphic, giving the conclusion of the lemma.
\end{proof}

\begin{Rk}\label{reduction}
In the context of Definition \ref{def}, suppose that the profinite group $G$ is finite, so that $G$ is a discrete space. 
Since the trivial subgroup $\{e\}$ is an open normal subgroup of $G$, \cite[Lemma 2.1.1]{Ribes} implies that for some 
$i_0 \in \{i\}$, $N_{i_0} = \{e\}$. Thus, $i_0$ is a terminal object of the directed poset $\{i\}$, so that if $\{Y_i\}_i$ is 
any inverse system of fibrant spectra indexed by $\{i\}$, then the natural map $\holim_i Y_i \to Y_{i_0}$ is a weak equivalence, by \cite[page 299: Example 4.1, (iii)]{Bousfield/Kan}. Therefore, given an $S[[G]]$-module $X = \holim_i X_i$, 
the natural map
\[\holim_i (X_i[(G/N_i)^l])_\mathtt{f} \xrightarrow{\simeq} (X_{i_0}[(G/N_{i_0})^l])_\mathtt{f}\] is a weak equivalence 
of spectra for each $l \geq 0$. By Lemma \ref{we}, this conclusion implies that
\[X_{hG} = d(\holim_i (X_i[(G/N_i)^\bullet])_\mathtt{f}) \xrightarrow{\simeq} d((X_{i_0}[(G/N_{i_0})^\bullet])_\mathtt{f})\]
is a weak equivalence. Also, there are weak equivalences
\[d((X_{i_0}[(G/N_{i_0})^\bullet])_\mathtt{f}) \xleftarrow{\simeq} d(X_{i_0}[G^\bullet]) \xleftarrow{\cong} (X_{i_0})_{h'G} 
\xleftarrow{\simeq} X_{h'G},\] where the leftmost equivalence is by Lemma \ref{we}, the isomorphism is 
due to Theorem \ref{model}, and the last equivalence is induced by the weak equivalence 
$X \xrightarrow{\simeq} X_{i_0}$, which is $G$-equivariant. Therefore, by a zigzag of weak equivalences, 
\[X_{hG} \simeq X_{h'G}\,,\] so that when $G$ is finite, $X_{hG}$ is equivalent to the usual homotopy orbit spectrum. 
\end{Rk}

\begin{Rk}\label{whenfinite}
Let $G$ be finite and let $X$ be any $G$-spectrum. Let $\{i\} = \{0\},$ a one-element set, and let 
$N_0 = 
\{e\},$ with $X_0=X_\mathtt{f}$. Then for $\{N_0\}_{0 \in \{0\}}$, the pair $(\{X_0\}_{0 \in \{0\}}, \holim_{\{0\}} X_0)$ is an $S[[G]]$-module, and the isomorphism 
$\holim_{\{0\}} X_0 \cong X_\mathtt{f}$ can be interpreted as saying that the $G$-spectrum $X$ can be realized as an $S[[G]]$-module. Also, we have 
\[(\holim_{\{0\}} X_0)_{hG} =
d(\holim_{\{0\}} (X_\mathtt{f}[G^\bullet])_\mathtt{f}) \xleftarrow{\cong} 
d((X_\mathtt{f}[G^\bullet])_\mathtt{f}) \xleftarrow{\simeq} d(X[G^\bullet]) \xleftarrow{\cong} X_{h'G},\] where the 
middle weak equivalence is by Lemma \ref{we}. 
Thus, there is a weak equivalence $X_{h'G} \xrightarrow{\simeq} 
(\holim_{\{0\}} X_0)_{hG}$, so that Definition \ref{def} recovers the classical 
definition of homotopy orbits.
\end{Rk}

If $K$ is a finite group and $Y$ is a spectrum with the trivial $K$-action, then the $K$-action 
on $K$ makes $Y[K]$ a $K$-spectrum and the augmentation 
$Y[K] \to Y$ induces a weak equivalence $(Y[K])_{h'K} \xrightarrow{\simeq} Y$. In 
Theorem \ref{Z[[G]]theorem} below, we 
show that there is a version of this result when the group at hand is a countably based 
profinite group.

\begin{Def}\label{Z[[G]]}
Let $G$ be a countably based profinite group, with 
\[N_0 \geq N_1 \geq \cdots \geq N_i \geq  
\cdots\] a chain of open normal subgroups of $G$, such that $G \cong \lim_{i \geq 0} G/N_i$. 
Also, let $Z$ be any spectrum. As an instance of a construction in \cite[page 375]{modular} and \cite[page 792]{GHMR}, 
we define the spectrum
\[Z[[G]] := \holim_{i \geq 0} (Z[G/N_i])_\mathtt{f}.\] We regard each $Z[G/N_i]$ as a $G$-spectrum by 
letting $G$ act trivially on $Z$ and in the usual way on $G/N_i$. It follows that with respect to 
the collection $\{N_i\}_{i \geq 0}$, the pair $(\{(Z[G/N_i])_\mathtt{f}\}_{i \geq 0}, \holim_{i \geq 0} (Z[G/N_i])_\mathtt{f})$ is a countably based $S[[G]]$-module. 
\end{Def} 

In the result below, we write ``$\,\simeq\,$" for a zigzag of weak equivalences whose exact form 
is specified in the proof. We note that the proof was 
aided by helpful comments from the referee of the revised versions of this paper. 

\begin{Thm}\label{Z[[G]]theorem}
If $G$ is a countably based profinite group and $Z$ is a spectrum, then there is an 
equivalence
\[(Z[[G]])_{hG} \simeq Z.\] 
\end{Thm}

\begin{proof} 
First, we recall the standard way that the weak equivalence $(Y[K])_{h'K} \xrightarrow{\simeq} Y$ 
mentioned above is justified. The 
augmentation $Y[K] \to Y$ extends to equip the simplicial spectrum $(Y[K])[K^\bullet]$ (see 
Definition \ref{xgbullet}) with an augmentation, giving a diagram of the form
\[\xymatrix@C.25in{
\cdots \ar@<.7ex>[r] \ar[r] \ar@<-.7ex>[r] & (Y[K])[K] \ar@<.5ex>[r] \ar@<-.5ex>[r] & (Y[K])[\ast] 
\ar[r] & Y.}\] The inclusion $Y \to Y[K] \xrightarrow{\cong} (Y[K])[\ast]$, where the first map 
sends $Y$ to the copy of $Y$ in the coproduct indexed by the identity element of $K$, is the 
``degree $-1$ map" of a contracting homotopy of the augmented simplicial spectrum $(Y[K])[K^\bullet] \to Y$, which implies that the induced composition \[(Y[K])_{h'K} \xrightarrow{\cong} 
d((Y[K])[K^\bullet]) \xrightarrow{\simeq} Y\] 
is a weak equivalence. In the preceding recollection, all the constructions are natural in $K$ and 
$Y$. 

By the above discussion, there is a tower
\[\Bigl\{(Z[G/N_i])[(G/N_i)^\bullet] \to Z\Bigr\}_{i \geq 0}\] of augmented simplicial 
spectra, with each augmented simplicial spectrum having a natural contracting homotopy. Here, 
$\{Z\}_{i \geq 0}$ is the constant tower on $Z$. Then 
by the functoriality of fibrant replacement and $\holim_{i}(-)$, there is the augmented simplicial spectrum 
\[\holim_{i}((Z[G/N_i])[(G/N_i)^\bullet])_\mathtt{f} \to 
\holim_{i} Z_\mathtt{f}\] and it has a contracting homotopy, which yields a 
weak equivalence
\[\mathcal{D} := d(\holim_{i}((Z[G/N_i])[(G/N_i)^\bullet])_\mathtt{f}) 
\xrightarrow{\simeq} \holim_{i} Z_\mathtt{f},\] whose source we write as $\mathcal{D}$ to 
simplify our notation. This weak equivalence fits into 
the zigzag of weak equivalences
\begin{equation}\label{mapI}\zig
\mathcal{D} \xrightarrow{\simeq} \holim_{i} Z_\mathtt{f} 
\xleftarrow{\simeq} Z_\mathtt{f} \xleftarrow{\simeq} Z,\end{equation} 
where the middle equivalence is 
by \cite[5.40]{Thomason} and the rightmost equivalence is the fibrant replacement map. 

For any $i, l \geq 0$, the fibrant replacement $Z[G/N_i] \to (Z[G/N_i])_\mathtt{f}$ induces 
a weak equivalence 
\[((Z[G/N_i])[(G/N_i)^l])_\mathtt{f} \xrightarrow{\simeq} ((Z[G/N_i])_\mathtt{f}[(G/N_i)^l])_\mathtt{f}\] between 
fibrant spectra, so that for each $l \geq 0$, there is a weak equivalence 
\[\holim_{i} ((Z[G/N_i])[(G/N_i)^l])_\mathtt{f} \xrightarrow{\simeq} 
\holim_{i} ((Z[G/N_i])_\mathtt{f}[(G/N_i)^l])_\mathtt{f}.\] 
Then it follows from Lemma \ref{we} that there is a weak equivalence
\[\underbrace{d(\holim_{i}((Z[G/N_i])[(G/N_i)^\bullet])_\mathtt{f})}_{\displaystyle{= \mathcal{D}}} \, \xrightarrow{\,\simeq\,} \,
\underbrace{d(\holim_{i}((Z[G/N_i])_\mathtt{f}[(G/N_i)^\bullet])_\mathtt{f})}_{\displaystyle{= (Z[[G]])_{hG}}}.
\] 
Putting this last weak equivalence together with (\ref{mapI}) gives the desired 
zigzag of 
weak equivalences.
\end{proof}

\end{section}

\section{Conditions for homotopy orbits to preserve weak equivalences}

We begin by recalling from Definition \ref{module} the concept ``morphism of $S[[G]]$-modules" and we define how to 
compose these morphisms. 

\begin{Def}\label{strong}
Let $G$ be any profinite group and let $\{N_i\}_i$ be a fixed cofinal subcollection of all the open normal subgroups of $G$. Let 
\[(\tau, \holim_i \tau_i) \: (\{X_i\}_i, \holim_i X_i) \to (\{Y_i\}_i, \holim_i Y_i)\] be a morphism of 
$S[[G]]$-modules. Thus, $\tau \: \{X_i\}_i \to \{Y_i\}_i$ is a natural transformation of diagrams of 
$G$-spectra, with component $\tau_i \: X_i \to Y_i$ for each $i$, and with $X = \holim_i X_i$ and $Y = \holim_i Y_i$, $\holim_i \tau_i$ is the 
induced $G$-equivariant map 
\[\holim_i \tau_i \: X = \holim_i X_i \to \holim_i Y_i = Y\] of spectra. Also, let 
\[(\tau', \holim_i \tau'_i) \: (\{Y_i\}_i, \holim_i Y_i) \to (\{Z_i\}_i, \holim_i Z_i)\] be a morphism of 
$S[[G]]$-modules. The natural transformation 
$\tau' \: \{Y_i\}_i \to \{Z_i\}_i$ has components $\tau'_i \: Y_i \to Z_i$. Then the composition 
$(\tau', \holim_i \tau'_i) \circ (\tau, \holim_i \tau_i)$ is defined to be the 
morphism 
\[(\{\tau'_i \circ \tau_i\}_i, \holim_i (\tau'_i \circ \tau_i)) \: (\{X_i\}_i, \holim_i X_i) \to (\{Z_i\}_i, \holim_i Z_i)\] of $S[[G]]$-modules. Here, $\{\tau'_i \circ \tau_i\}_i \: \{X_i\}_i \to \{Z_i\}_i$ is a natural 
transformation of diagrams of $G$-spectra, with components $\tau'_i \circ \tau_i \: X_i \to Z_i$. 
\end{Def}

\begin{Def}
In Definition \ref{strong}, for each $i$ and for every $l \geq 0$, $\tau_i$ 
induces a map \[(X_i[(G/N_i)^l])_\mathtt{f} \rightarrow (Y_i[(G/N_i)^l])_\mathtt{f},\] so that there is the 
induced map 
\[\holim_i (X_i[(G/N_i)^\bullet])_\mathtt{f} \rightarrow \holim_i (Y_i[(G/N_i)^\bullet])_\mathtt{f}\] of simplicial spectra. 
By applying $d(-)$ to this map, we see that $\tau$ induces 
a map $X_{hG} \to Y_{hG}$ that in a slight abuse of notation, we denote by \[\tau_{hG} \: X_{hG} \to Y_{hG}.\] 
\end{Def} 

\begin{Rk}\label{omitagain}
Given a profinite group $G$, in the expression ``let
\[(\tau, \holim_i \tau_i) \: (\{X_i\}_i, \holim_i X_i) \to (\{Y_i\}_i, \holim_i Y_i)\] be a morphism 
of $S[[G]]$-modules" and its natural variants, we omit mentioning the fixed family $\{N_i\}_i$ that is a part of the picture. 
\end{Rk}

In Definition \ref{strong}, if $\tau_i$ is a weak equivalence of spectra for each $i$, then the map 
$\holim_i \tau_i$ is a weak equivalence of spectra. The following result shows that in this case, 
the map $X_{hG} \to Y_{hG}$ is also a weak equivalence.

\begin{Thm}\label{induces}
Let $G$ be a profinite group and let 
\[(\tau, \holim_i \tau_i) \: (\{X_i\}_i, \holim_i X_i) \to (\{Y_i\}_i, \holim_i Y_i)\] be a morphism 
of $S[[G]]$-modules, with $X = \holim_i X_i$ and $Y = \holim_i Y_i$. 
If for each $i$, $\tau_i \: X_i \to Y_i$ is a weak equivalence of spectra, then 
the map $\tau_{hG} \: X_{hG} \rightarrow Y_{hG}$ is a weak equivalence of spectra.
\end{Thm}
\begin{proof}
Let $l \geq 0$. For each $i$, since $X_i \rightarrow Y_i$ is a weak equivalence, the map 
\[(X_i[(G/N_i)^l])_\mathtt{f} \rightarrow (Y_i[(G/N_i)^l])_\mathtt{f}\] is a weak equivalence 
between fibrant spectra. Thus,  
\[\holim_i (X_i[(G/N_i)^l])_\mathtt{f} \rightarrow \holim_i (Y_i[(G/N_i)^l])_\mathtt{f}\] is a weak equivalence, so that by 
Lemma \ref{we}, $X_{hG} \rightarrow Y_{hG}$ is a weak equivalence.
\end{proof} 

Now we give a generalization of Theorem \ref{induces}. We 
separated out Theorem \ref{induces} and gave its proof because the result describes a 
natural situation and its justification is brief. 

\begin{Thm}\label{generalized}
Let $G$ be any profinite group and let \[(\tau, \holim_i \tau_i) \:  (\{X_i\}_i, \holim_i X_i) \to (\{Y_i\}_i, \holim_i Y_i)\] be a 
morphism of $S[[G]]$-modules, with $X = \holim_i X_i$ and $Y = \holim_i Y_i$. If the 
map $\holim_i \tau_i$ is a weak equivalence of spectra, then the map $\tau_{hG} \: 
X_{hG} \to Y_{hG}$ 
is a weak equivalence of spectra. 
\end{Thm}
\begin{proof}
Let $l \geq 0$. 
As in the 
proof of Theorem \ref{induces}, Lemma \ref{we} implies that we only need to show that the induced map 
\[\holim_i (X_i[(G/N_i)^l])_\mathtt{f} \to \holim_i (Y_i[(G/N_i)^l])_\mathtt{f}\] of 
spectra is a weak equivalence. 

By relabeling, we can write the directed poset $\{i\}$ as $\{k\}$ and also as $\{i'\}$. Let 
\[\Delta := \{(i, i) \in \{k\} \times \{i'\} \mid i \in \{k\}\}\] be the diagonal of $\{i\} \times \{i\}$. 
Notice that there is a canonical map
\[\holim_i (X_i[(G/N_i)^l])_\mathtt{f} = \holim_{(i,i) \in \Delta} (X_i[(G/N_i)^l])_\mathtt{f} 
\xleftarrow{\,\simeq\,} \holim_{(i,j) \in \{k\} \times \{i'\}} (X_i[(G/N_j)^l])_\mathtt{f},\] and it 
is a weak equivalence because in the terminology of \cite[Chapter XI, \S 9.1]{Bousfield/Kan}, 
the inclusion functor 
$\Delta^\mathrm{op} \to (\{k\} \times \{i'\})^\mathrm{op}$ is left cofinal. This left cofinality 
is easy to see by applying \cite[Chapter XI, Proposition 9.3]{Bousfield/Kan}. 
Also, there is 
the natural isomorphism
\[\holim_{j \in \{i'\}} \holim_{i \in \{k\}} (X_i[(G/N_j)^l])_\mathtt{f} \xrightarrow{\cong} 
\holim_{(i,j) \in \{k\} \times \{i'\}} (X_i[(G/N_j)^l])_\mathtt{f}.\] To simplify our notation, 
we write the expression $\displaystyle{\holim_{j \in \{i'\}} \holim_{i \in \{k\}} (X_i[(G/N_j)^l])_\mathtt{f}}$
as $\displaystyle{\holim_{j} \holim_{i} (X_i[(G/N_j)^l])_\mathtt{f}}$. 
The preceding remarks go 
through for the $S[[G]]$-module $(\{Y_i\}_i, \holim_i Y_i)$ as well, and hence, we obtain the commutative square 
\[\xymatrix{
\displaystyle{\holim_{j} \holim_{i} (X_i[(G/N_j)^l])_\mathtt{f}} \ar[r]^-\simeq \ar[d] & \displaystyle{\holim_i (X_i[(G/N_i)^l])_\mathtt{f}} \ar[d] \\ 
\displaystyle{\holim_{j} \holim_{i} (Y_i[(G/N_j)^l])_\mathtt{f}} \ar[r]^-\simeq & \displaystyle{\holim_i (Y_i[(G/N_i)^l])_\mathtt{f}}, 
}\]
in which the horizontal arrows are weak equivalences. Thus, to complete the proof, it 
suffices to show that the left vertical arrow is a weak equivalence.

Let $i$ and $j$ be arbitrary and notice that $(G/N_j)^l$ is finite. Then there is the commutative diagram
\[\xymatrix{
(X_i[(G/N_{j})^l])_\mathtt{f} \ar[d] & ({\bigvee}_{(G/N_{j})^l} X_i)_\mathtt{f} \ar[l]_-\simeq 
\ar[r]^-\simeq \ar[d] & (\prod_{(G/N_{j})^l} X_i)_\mathtt{f} \ar[d] \\
(Y_i[(G/N_{j})^l])_\mathtt{f} & {(\bigvee_{(G/N_{j})^l} Y_i)_\mathtt{f}} \ar[l]_-\simeq 
\ar[r]^-\simeq & {(\prod_{(G/N_{j})^l} Y_i)_\mathtt{f}}
}\] built from fibrant replacement of canonical maps, in which every horizontal arrow is a weak equivalence between fibrant spectra. By applying $\holim_i (-)$, it 
follows that for each $j$, there is the commutative diagram
\[\xymatrix{
\displaystyle{\holim_i (X_i[(G/N_{j})^l])_\mathtt{f}} \ar[d] & \displaystyle{\holim_i} \textstyle{({\bigvee}_{(G/N_{j})^l} X_i)_\mathtt{f}} \ar[l]_-\simeq 
\ar[r]^-\simeq \ar[d] & \displaystyle{\holim_i} \textstyle{(\prod_{(G/N_{j})^l} X_i)_\mathtt{f}} \ar[d] \\
\displaystyle{\holim_i (Y_i[(G/N_{j})^l])_\mathtt{f}} & \displaystyle{\holim_i} \textstyle{{(\bigvee_{(G/N_{j})^l} Y_i)_\mathtt{f}}} \ar[l]_-\simeq 
\ar[r]^-\simeq & \displaystyle{\holim_i} \textstyle{{(\prod_{(G/N_{j})^l} Y_i)_\mathtt{f}}}, 
}\] with each horizontal arrow a weak equivalence. Notice that the left vertical arrow is a 
weak equivalence if the rightmost vertical arrow is a weak equivalence. 

We continue to let $j$ be arbitrary. For each $i$, $X_i$ and $Y_i$ are fibrant spectra, so that $\prod_{(G/N_j)^l} X_i$ and $\prod_{(G/N_j)^l} Y_i$ are too. Then by using canonical maps and 
fibrant replacement, there is the 
commutative diagram
\[\xymatrix{
\displaystyle{\holim_i} (\textstyle{\prod_{(G/N_{j})^l} X_i})_\mathtt{f} \ar[d] & \displaystyle{\holim_i} \,\textstyle{\prod_{(G/N_{j})^l} X_i}
\ar[l]_-\simeq \ar[r]^-\cong \ar[d] & \prod_{(G/N_{j})^l} \displaystyle{\holim_i X_i}  \ar[d]\\
\displaystyle{\holim_i}\textstyle{(\prod_{(G/N_{j})^l} Y_i)_\mathtt{f}} & \displaystyle{\holim_i} \,\textstyle{\prod_{(G/N_{j})^l} Y_i}
\ar[l]_-\simeq \ar[r]^-\cong & \prod_{(G/N_{j})^l} \displaystyle{\holim_i Y_i}\,
}\] in which the horizontal arrows in each row are, on the left, a weak equivalence and, 
on the right, an isomorphism. Since $\holim_i X_i \to \holim_i Y_i$ is a weak equivalence, 
the rightmost vertical arrow in our last diagram above is a weak equivalence, giving that the 
leftmost vertical arrow in this diagram is one too, and hence, for each $j$, the map 
\[\holim_i (X_i[(G/N_{j})^l])_\mathtt{f} \to \holim_i (Y_i[(G/N_{j})^l])_\mathtt{f}\] is a 
weak equivalence between fibrant spectra. It follows that the map
\[\holim_j \holim_i (X_i[(G/N_{j})^l])_\mathtt{f} \to \holim_j \holim_i (Y_i[(G/N_{j})^l])_\mathtt{f}\]  
is a weak equivalence. 
\end{proof}

The following result describes three situations in which the hypotheses of Theorem 
\ref{generalized} are satisfied. 

\begin{Cor}\label{3situations}
Let $G$ be a profinite group, with \[(\tau, \holim_i \tau_i) \: (\{X_i\}_i, \holim_i X_i)  \to (\{Y_i\}_i, \holim_i Y_i)\] a morphism of $S[[G]]$-modules, where $X = \holim_i X_i$ and $Y = \holim_i Y_i$. 
Also, suppose that
\[\lim_i \pi_t(\tau_i) \: \lim_i \pi_t(X_i) \to \lim_i \pi_t(Y_i)\] is a bijection for every integer $t$. Consider the conditions
\begin{itemize}
\item[{(i)}] 
for every integer $t$, the 
inverse systems $\{\pi_t(X_i)\}_i$ and $\{\pi_t(Y_i)\}_i$ consist of compact Hausdorff 
abelian groups and continuous homomorphisms, 
\item[{(ii)}] 
$G$ is countably based and the fixed family $\{N_i\}_i$ is equal to a chain $\{N_i\}_{i \geq 0}$, as in 
$\mathrm{(\ref{chainofopens})}$,
\item[{(iii)}]
for every $t$, the canonical map $\lim^1_i \pi_t(X_i) \to \lim^1_i \pi_t(Y_i)$ is an isomorphism, and
\item[{(iv)}]  
for every $t$ and all $i$, $\pi_t(X_i)$ and 
$\pi_t(Y_i)$ are finitely generated abelian groups.   
\end{itemize} 
If $\mathrm{(i)}$, or the pair $\mathrm{(ii)}$ and $\mathrm{(iii)}$, or the pair $\mathrm{(iii)}$ and $\mathrm{(iv)}$ is satisfied, then both $\holim_i \tau_i$ and $\tau_{hG} \: X_{hG} \to Y_{hG}$ are weak equivalences of spectra. 
\end{Cor}

\begin{proof}
There is the morphism
\[\xymatrix{
E_2^{s,t}(\{X_i\}_i) \ar@{=>}[r] \ar[d]_-{\lambda^{s,t}} & \pi_{t-s}(\displaystyle{\holim_i X_i}) \ar[d]\\
E_2^{s,t}(\{Y_i\}_i) \ar@{=>}[r] & \pi_{t-s}(\displaystyle{\holim_i Y_i})
}\] of homotopy spectral sequences, where the map $\lambda^{s,t}$ is the induced homomorphism 
$\lim^s_i \pi_t(X_i) \to \lim^s_i \pi_t(Y_i)$ of $E_2$-terms. It follows that when the bigraded map $\lambda^{\ast,\ast}$ is 
an isomorphism, $\holim_i \tau_i$ is a weak equivalence. 

To obtain the desired conclusion, we only have to show that the graded map 
$\lambda^{s,\ast}$ is an isomorphism in each of the following three distinct cases:
\begin{itemize}
\item 
for $s \geq 1$, when condition (i) holds;
\item 
for $s \geq 2$, when (ii) and (iii) hold;
\item
for $s \geq 2$, when (iii) and (iv) hold. 
\end{itemize} 
Thus, the desired conclusion is an immediate consequence of the following. By 
\cite[Theorem 2]{Yosimura} and \cite{Jensenbook}, when (i) holds, 
\[\lims_i \pi_t(X_i) = \lims_i \pi_t(Y_i) = 0, \ \ \ \text{for} \ s \geq 1.\] When (ii) or (iv) is satisfied, 
for each $s \geq 2$, $\lim^s_i \pi_t(X_i) = \lim^s_i \pi_t(Y_i) = 0$. In the case of (ii), this 
conclusion is well-known, and in the case of (iv), it holds by \cite[page 65]{Jensenbook}.
\end{proof}
 
\begin{section}{The homotopy orbit spectral sequence for an $S[[G]]$-module}
\label{SS}

Our work in this section begins with recalling the classical homotopy orbit 
spectral sequence for $X_{h'G}$, when $G$ is finite, and doing some related preparatory work. If 
$A$ is an abelian group and $K$ is a group, then we use the notation 
\[A[K] := \textstyle{\bigoplus_K A}.\]

Let $G$ be 
a finite group and let $X$ be a $G$-spectrum.
By Theorem \ref{dss}, there is a spectral sequence 
\[E_2^{p,q} 
\Longrightarrow \pi_{p+q}(d(X[G^\bullet])) = \pi_{p+q}(X_{h'G}),\] where
\begin{equation}\zig\label{formula}
E_2^{p,q} = H_p(\pi_q(X[G^\ast])) \cong H_p(G,\pi_q(X)),
\end{equation} the 
$p$th group homology of $G$, with coefficients in $\pi_q(X)$ (see, for example, 
\cite[(7.9)]{Jardine}). We recall from Section \ref{two} that for each $k \geq 0$, there is a simplicial 
pointed simplicial set $X_k[G^\bullet]$ that satisfies the isomorphism
\begin{equation}\label{isowithreplace}\zig
(X[G^\bullet])_k = X_k[G^\bullet] \cong \textstyle{\coprod_\ast (G \to X_k)}\end{equation} 
where the last 
expression is the simplicial 
replacement in $\mathcal{S}_\ast$, the category of pointed simplicial sets. 

It is helpful to recall from 
\cite[Chapter XII, \S 5.7]{Bousfield/Kan} that if $\mathcal{C}$ is any small category and 
$\underline{Y}$ is a diagram $\mathcal{C} \to \mathcal{S}_\ast$, then given (to quote \cite[page 339]{Bousfield/Kan}) a `reduced generalized homology theory $\tilde{h}_\ast$ which ``comes from a spectrum,"' there is a 
natural isomorphism
\begin{equation}\label{B&Kimplicit}\zig
\tilde{h}_t(\textstyle{\coprod_\ast \underline{Y}}) \cong \textstyle{\coprod_\ast \tilde{h}_t(\underline{Y})}, \ \ \ t \in \mathbb{Z}
\end{equation} 
of simplicial abelian groups, 
where $\coprod_\ast \underline{Y}$ is the usual simplicial replacement of 
$\underline{Y}$ in $\mathcal{S}_\ast$ and $\coprod_\ast \tilde{h}_t(\underline{Y})$ is the 
simplicial replacement in $\mathbf{Ab}$, the category of abelian groups, 
of the diagram 
$\tilde{h}_t(\underline{Y}) \: \mathcal{C} \to \mathbf{Ab}$ (see \cite[Chapter XII, \S 5.4, \S 5.5]{Bousfield/Kan}). The simplicial replacement $\coprod_\ast \tilde{h}_t(\underline{Y})$ is defined as in 
$\mathcal{S}_\ast$, but with the coproduct $\bigoplus$ in $\mathbf{Ab}$ replacing the coproduct $\bigvee$ in $\mathcal{S}_\ast$. For more detail about simplicial replacement in categories other than just $\mathcal{S}_\ast$, see for example, \cite[pages 211-212]{radditive} and \cite[2.5, 2.11]{Beatriz}. 

Let $q \in \mathbb{Z}$ and recall that if $Z$ is any spectrum, $\pi_q(Z) = \displaystyle{\colim_{k \geq 0,\,q+k \geq 2} \pi_{q+k}(Z_k)}$. Now we consider $X[G^\bullet]$ from above further. There are natural isomorphisms
\begin{align*}
\pi_q(X[G^\bullet]) & = \colim_{k \geq 0,\,q+k \geq 2} \pi_{q+k}((X[G^\bullet])_k) 
\cong \colim_{k \geq 0,\,q+k \geq 2} \textstyle{\pi_{q+k}(\coprod_\ast (G \to X_k))}\\
& \cong \colim_{k \geq 0,\,q+k \geq 2} (\textstyle{\coprod_\ast (G \to \pi_{q+k}(X_k)))}\end{align*} 
of simplicial abelian groups, where $\coprod_\ast (G \to \pi_{q+k}(X_k))$ is the simplicial replacement 
in $\mathbf{Ab}$ of the diagram $G \to \pi_{q+k}(X_k)$ defined by the induced action of $G$ on $\pi_{q+k}(X_k)$, the first isomorphism is by (\ref{isowithreplace}), and the 
second isomorphism is by (\ref{B&Kimplicit}). In each dimension $l \geq 0$ of the above 
$\coprod_\ast (G \to \pi_{q+k}(X_k))$, 
\[\textstyle{\coprod_l (G \to \pi_{q+k}(X_k)) = \bigoplus_{G^l} \pi_{q+k}(X_k)},\] so that 
there is a natural isomorphism
\begin{align*}
\colim_{k \geq 0,\,q+k \geq 2} (\textstyle{\coprod_l (G \to \pi_{q+k}(X_k)))} & \cong 
\textstyle{\bigoplus_{G^l}} \, \displaystyle{\colim_{k \geq 0,\,q+k \geq 2} \pi_{q+k}(X_k)} = \textstyle{\bigoplus_{G^l} \pi_q(X)} \\ &= \textstyle{\coprod_l (G \to \pi_q(X))},
\end{align*} 
where the last expression is dimension $l$ of $\coprod_\ast (G \to \pi_q(X))$, the 
simplicial replacement in $\mathbf{Ab}$ of the diagram $G \to \pi_q(X)$, defined by the induced 
action 
of $G$ on $\pi_q(X)$. By putting together the 
various natural isomorphisms above, we obtain the isomorphisms 
\begin{equation}\label{isoInsAb}\zig
\pi_q(X[G^\bullet]) \cong \colim_{k \geq 0,\,q+k\geq 2} (\textstyle{\coprod_\ast (G \to \pi_{q+k}(X_k)))} 
\cong \textstyle{\coprod_\ast (G \to \pi_q(X))}\end{equation} of simplicial abelian groups. 
\begin{Def}\label{pixgbullet}
Let $G$ be a finite group and suppose that $X$ is a $G$-spectrum. By (\ref{isoInsAb}), for every integer $q$ and each $l \geq 0$, there is a natural isomorphism 
\[\pi_q(X[G^l]) \cong \pi_q(X)[G^l] = \textstyle{\bigoplus}_{G^l} \pi_q(X) = \coprod_l (G \to \pi_q(X)),\] and so we define the simplicial abelian group
\[\textstyle{\pi_q(X)[G^\bullet] := \coprod_\ast (G \to \pi_q(X)).}\] Following Definition \ref{Moore}, 
the Moore complex of 
$\pi_q(X)[G^\bullet]$ is denoted by $\pi_q(X)[G^\ast]$, and the isomorphisms of simplicial abelian groups in (\ref{isoInsAb}) imply that there is a natural isomorphism
\begin{equation}\label{helpfuliso}\zig
\pi_q(X[G^\ast]) \cong \pi_q(X)[G^\ast]\end{equation} of chain complexes. 
\end{Def}

In the context of Definition \ref{pixgbullet}, notice that by (\ref{formula}), there is the isomorphism 
\begin{equation}\zig\label{homology} 
H_p(G,\pi_q(X)) \cong H_p(\pi_q(X)[G^\ast]), \ \ \ \text{for} \ p \geq 0, \ q \in \mathbb{Z}.
\end{equation}

Now we are ready to study the homotopy orbit spectral sequence when $G$ is profinite. 
The following result is an immediate application of Definition \ref{def} and 
Theorem \ref{dss}. 

\begin{Thm}\label{generalh-oSS}
Let $G$ be any profinite group and let 
$X = \holim_i X_i$ be an $S[[G]]$-module. Then there is a 
spectral sequence having the form
\[E_2^{p,q} = H_p(\pi_q(\holim_i (X_i[(G/N_i)^\ast])_\mathtt{f})) \Longrightarrow \pi_{p+q}(X_{hG}),\] where $E_2^{p,q}$ is the $p$th homology of the specified Moore complex. 
\end{Thm}

\begin{Rk}
By Remark \ref{bar2}, the $E_2$-term in Theorem \ref{generalh-oSS} can be written 
more succinctly as
\[E_2^{p,q} = \pi_p\pi_q\mspace{-3mu}\left(B^\mathrm{cts}_\bullet(G,X)\right),\] the $p$th homotopy group of the stated simplicial abelian group.
\end{Rk}

Comparing the $E_2$-term of the homotopy orbit spectral sequence 
in Theorem \ref{generalh-oSS} with the $E_2$-term in (\ref{formula}) 
motivates one to try to 
identify situations in which the former $E_2$-term can 
be described in a more meaningful homological way. To help with this, we have the following 
result.  

\begin{Thm}\label{compactHaus}
Let $G$ be any profinite group and let $X = \holim_i X_i$ be an $S[[G]]$-module. Given an 
integer $t$ and $l \geq 0$, let 
\[\psi(t,l) \: \pi_t(\holim_i(X_i[(G/N_i)^l])_\mathtt{f}) \to
\lim_i \pi_{t}(X_i)[(G/N_i)^l]\] be the canonical map whose source is 
the $t^\text{th}$ homotopy group of the $l$-simplices of the simplicial spectrum 
$\holim_i (X_i[(G/N_i)^\bullet])_\mathtt{f}$. If for every $t$, the 
inverse system $\{\pi_t(X_i)\}_i$ consists of compact Hausdorff abelian groups and 
continuous homomorphisms, then for every $t$ and all $l \geq 0$, the map $\psi(t,l)$ 
is an isomorphism of abelian groups.
\end{Thm}

\begin{proof}
Let $l \geq 0$. The canonical map is given by the universal property of the limit and the fact that for each $i$, 
since $(G/N_i)^l$ is finite, there is an isomorphism
\[\pi_\ast((X_i[(G/N_i)^l])_\mathtt{f}) \cong \pi_\ast(X_i)[(G/N_i)^l].\] There is the homotopy spectral sequence 
\[E^{s,t}_2  \Longrightarrow \pi_{t-s}(\holim_i (X_i[(G/N_i)^l])_\mathtt{f}),\] with
\[
E_2^{s,t} = \lims_i \pi_t((X_i[(G/N_i)^l])_\mathtt{f}) \cong \lims_i \pi_t(X_i)[(G/N_i)^l].\]  

Let $t$ be any integer. Since each $(G/N_i)^l$ is finite and the category $\mathcal{CHA}$ of compact Hausdorff abelian groups is abelian, the direct sums $\pi_t(X_i)[(G/N_i)^l]$ are coproducts in $\mathcal{CHA}$, so that the inverse system 
$\{\pi_t(X_i)[(G/N_i)^l]\}_i$ consists of compact Hausdorff abelian groups. Also, all homomorphisms
\[\pi_t(X_i)[(G/N_i)^l] \to \pi_t(X_{i'})[(G/N_{i'})^l]\] in the inverse system are given by the 
universal property of the coproduct in $\mathcal{CHA}$, and hence, these homomorphisms 
are continuous. For example, if $g \in (G/N_i)^l$ and $g'$ is its 
image under the natural map $(G/N_i)^l \to (G/N_{i'})^l$, there is the commutative diagram 
\[\xymatrix{
\bigoplus_{(G/N_i)^l} \pi_t(X_i) \ar[r] & \bigoplus_{(G/N_{i'})^l} \pi_t(X_{i'})\\
\pi_t(X_i) \ar[u] \ar[r] & \pi_t(X_{i'}), \ar[u]
}\] in which the vertical maps are inclusions of the copies of $\pi_t(X_i)$ and $\pi_t(X_{i'})$ indexed by $g$ and $g'$, respectively, and the upper horizontal arrow is given by the universal property. This diagram is in $\mathcal{CHA}$ and, in particular, consists of continuous homomorphisms. We have shown that the inverse system 
$\{\pi_t(X_i)[(G/N_i)^l]\}_i$ lives in $\mathcal{CHA}$, so that by \cite[Theorem 2]{Yosimura} and \cite{Jensenbook}, 
\[\lims_i \pi_t(X_i)[(G/N_i)^l] = 0, \ \ \ s \geq 1,\] and hence, the above spectral sequence
collapses, giving that for every $t \in \mathbb{Z}$, the desired map
is an isomorphism.  
\end{proof}

Now we use Theorem \ref{compactHaus} to describe a situation in which the $E_2$-term of 
the homotopy orbit spectral sequence of Theorem \ref{generalh-oSS} can be identified in a homologically interesting way.  

\begin{Thm}\label{corollary:whenfinite}
Let $G$ be a profinite group and suppose that $X = \holim_i X_i$ is an $S[[G]]$-module. Suppose that for every integer $t$, $\{\pi_t(X_i)\}_i$ is an inverse system of compact Hausdorff abelian groups 
and continuous homomorphisms. Also, suppose that for all $t$ and each $i$, the induced action of the discrete group $G/N_i$ on the compact Hausdorff space 
$\pi_t(X_i)$ is continuous. Then there is a homotopy orbit spectral sequence of the 
form
\[E_2^{p,q} \cong \lim_i H_p(G/N_i, \pi_q(X_i)) \Longrightarrow \pi_{p+q}(X_{hG}).\] 
\end{Thm}

\begin{proof}
Let $p \geq 0$ and let $q$ be any integer. By Theorem \ref{generalh-oSS}, we only need to show that there is an isomorphism 
\[H_p(\pi_q(\holim_i (X_i[(G/N_i)^\ast])_\mathtt{f})) \cong \lim_i H_p(G/N_i, \pi_q(X_i)).\] 
There are isomorphisms
\[\pi_q(\holim_i (X_i[(G/N_i)^\ast])_\mathtt{f}) \cong \lim_i \pi_q(X_i[(G/N_i)^\ast]) \cong 
\lim_i (\pi_q(X_i)[(G/N_i)^\ast])\] of chain complexes, where the first step is by 
Theorem \ref{compactHaus} and its proof, and the second step applies (\ref{helpfuliso}). Thus, 
there is an isomorphism 
\[H_p(\pi_q(\holim_i (X_i[(G/N_i)^\ast])_\mathtt{f})) \cong H_p\Bigl[\lim_i (\pi_q(X_i)[(G/N_i)^\ast])\Bigr]. 
\] 

As explained in the proof of Theorem \ref{compactHaus}, for each $l \geq 0$ and every $i$, the abelian group $\pi_q(X_i)[(G/N_i)^l]$ is compact Hausdorff. 
Then our next step is to show that the inverse system
\[\{\pi_q(X_i)[(G/N_i)^\ast]\}_i\] of chain complexes lives in the category $\mathcal{CHA}$ of compact Hausdorff abelian groups. The proof of Theorem \ref{compactHaus} verified that for each $l \geq 0$, $\{\pi_q(X_i)[(G/N_i)^l]\}_i$ is an inverse system in $\mathcal{CHA}$. Thus, to complete the next step, it suffices 
to show that for each $i$, all the boundary homomorphisms of the chain complex $\pi_q(X_i)[(G/N_i)^\ast]$ are in 
$\mathcal{CHA}$.
Let $i$ be arbitrary. 
The chain complex $\pi_q(X_i)[(G/N_i)^\ast]$ is the Moore complex of 
the simplicial abelian group $\pi_q(X_i)[(G/N_i)^\bullet]$, 
and 
since $\mathcal{CHA}$ is abelian, the boundary homomorphisms of $\pi_q(X_i)[(G/N_i)^\ast]$ are in 
$\mathcal{CHA}$, if for every 
$l \geq 0$ and all $j$ such that $0 \leq j \leq l+1$, the face maps 
\[d_j \: \pi_q(X_i)[(G/N_i)^{l+1}] \to \pi_q(X_i)[(G/N_i)^l]\] of $\pi_q(X_i)[(G/N_i)^\bullet]$ 
are continuous. 

We continue to let $i$, $l$, and $j$ be as above. Recall that by Definition \ref{pixgbullet}, the simplicial abelian group $\pi_q(X_i)[(G/N_i)^\bullet]$ is $\coprod_\ast (G/N_i \to \pi_q(X_i))$, a 
simplicial replacement in $\mathbf{Ab}$, instead of in $\mathcal{S}_\ast$. To verify that each $d_j$ is continuous, we show that this simplicial replacement can be enriched: the whole construction can be carried out in the category $\mathcal{CHA}$. As in the proof of Theorem \ref{compactHaus}, the finite direct sums $\pi_q(X_i)[(G/N_i)^{l+1}]$ and 
$\pi_q(X_i)[(G/N_i)^l]$ are coproducts in $\mathcal{CHA}$: for example, $\pi_q(X_i)[(G/N_i)^{l+1}] = 
\coprod_{(G/N_i)^{l+1}} \pi_q(X_i)$ in $\mathcal{CHA}$ and the inclusion $\pi_q(X_i) \to \pi_q(X_i)[(G/N_i)^{l+1}]$, from the copy of 
$\pi_q(X_i)$ indexed by any $g \in (G/N_i)^{l+1}$ into the coproduct, is continuous. Also, since 
the discrete group $G/N_i$ acts continuously on $\pi_q(X_i)$, the action map
\[G/N_i \times \pi_q(X_i) \to \pi_q(X_i), \ \ \ (\gamma, m) \mapsto \gamma \cdot m\] is continuous. Hence, given $\gamma' \in G/N_i$, the induced map $\bar{\gamma'} \: \pi_q(X_i) \to \pi_q(X_i)$ that 
is defined to be the composition
\[\pi_q(X_i) \to G/N_i \times \pi_q(X_i) \to \pi_q(X_i), \ \ \ m \mapsto (\gamma', m) \mapsto \gamma' \cdot m,\] is a composition of continuous functions and is thereby continuous. Thus, the diagram $G/N_i \to \pi_q(X_i)$, which thus far has been described as landing in $\mathbf{Ab}$, can be 
regarded as having target category equal to $\mathcal{CHA}$.

With the preceding remarks in hand and by repeatedly using the universal property of finite coproducts in $\mathcal{CHA}$, it is now straightforward to go through the definitions in Section \ref{two} 
of the face and degeneracy maps to see that for each $i$, there is a simplicial replacement 
\[\textstyle{\coprod_\ast^{\scriptscriptstyle{\mathcal{CHA}}}} (G/N_i \to \pi_q(X_i))\] in the category $\mathcal{CHA}$ 
of the diagram $G/N_i \to \pi_q(X_i)$, with 
all face and degeneracy maps continuous homomorphisms. Also, if we let $\mathbb{U} \: \mathtt{s}(\mathcal{CHA}) \to \mathtt{s}(\mathbf{Ab})$ be the forgetful functor from the category of simplicial objects in $\mathcal{CHA}$ to simplicial abelian groups, then we see that there are identities
\[\mathbb{U}\mspace{-3mu}\left(\textstyle{\coprod_\ast^{\scriptscriptstyle{\mathcal{CHA}}}} (G/N_i \to \pi_q(X_i))\right) = \textstyle{\coprod_\ast} 
(G/N_i \to \pi_q(X_i)) = \pi_q(X_i)[(G/N_i)^\bullet]\] of simplicial abelian groups. Thus, each 
diagram $\pi_q(X_i)[(G/N_i)^\bullet]$ belongs to $\mathcal{CHA}$, and we can conclude that the inverse system $\{\pi_q(X_i)[(G/N_i)^\ast]\}_i$ of chain 
complexes lives in $\mathcal{CHA}$. 

By, for example, \cite[Proposition 4]{Yosimura} and \cite[proof of Proposition 2.2]{RognesSverre}, 
the functor $\lim_i (-)$ applied to $\{i\}$-indexed inverse systems in $\mathcal{CHA}$ is exact, so that 
forming $\lim_i (-)$ commutes with 
taking homology in this category. It follows that 
there is an isomorphism 
\[H_p\Bigl[\lim_i (\pi_q(X_i)[(G/N_i)^\ast])\Bigr] 
\cong \lim_i H_p(\pi_q(X_i)[(G/N_i)^\ast]).\] 
Therefore, we have
\begin{align*}
H_p(\pi_q(\holim_i (X_i[(G/N_i)^\ast])_\mathtt{f})) & \cong \lim_i H_p(\pi_q(X_i)[(G/N_i)^\ast])\\
& \cong \lim_i H_p(G/N_i, \pi_q(X_i)),\end{align*} where 
the second isomorphism is by (\ref{homology}).  
\end{proof}

Suppose that the hypotheses of Theorem \ref{corollary:whenfinite} are satisfied. Then for every integer $q$, the canonical $G$-equivariant map $\pi_q(X) \to \lim_i \pi_q(X_i)$ is an isomorphism of abelian groups, by the $l=0$ case of Theorem \ref{compactHaus}. Now suppose further 
that $q$ is a fixed integer such that the inverse system $\{\pi_q(X_i)\}_i$ consists of profinite groups. 
For each $i$, 
since the $G$-action on $X_i$ factors through $G/N_i$, the induced action of $G$ on $\pi_q(X_i)$ is given by the 
composition
\[G \times \pi_q(X_i) \xrightarrow{\pi \times \mathrm{id}} G/N_i \times \pi_q(X_i) \to \pi_q(X_i),\] where $\pi$ is the canonical projection, $\mathrm{id}$ is the identity map, and the last map in the composition is given by the continuous action of $G/N_i$ on $\pi_q(X_i)$. This composition is continuous, so that $\pi_q(X_i)$ is a profinite $G$-module. Thus, for every $i$, $\pi_q(X_i)$ is a profinite $\widehat{\mathbb{Z}}[[G]]$-module, by \cite[Proposition 5.3.6, (c)]{Ribes}, and hence, 
$\lim_i \pi_q(X_i)$ is a profinite $\widehat{\mathbb{Z}}[[G]]$-module. 

The above discussion shows that 
if the hypotheses of Theorem 
\ref{corollary:whenfinite} are satisfied and, for some integer $q$, $\pi_q(X_i)$ is a profinite group for every $i$, then we have the following:
\begin{itemize}
\item
$\lim_i \pi_q(X_i)$ is a profinite $\widehat{\mathbb{Z}}[[G]]$-module; 
\item 
we can identify the $G$-module 
$\pi_q(X)$ with $\lim_i \pi_q(X_i)$; 
\item 
under the preceding identification, $\pi_q(X)$ is a profinite $\widehat{\mathbb{Z}}[[G]]$-module, so that the continuous homology groups $H^c_\ast(G, \pi_q(X))$ are well-defined; and 
\item
by \cite[Proposition 6.5.7]{Ribes}, there is an isomorphism
\begin{equation}\label{isowithHcts}\zig
\lim_i H_\ast(G/N_i, \pi_q(X_i)) \cong H^c_\ast(G, \pi_q(X)).\end{equation}
\end{itemize}
These observations yield Corollary \ref{profinitecase} below. Before stating this result, we provide some additional motivation for (\ref{isowithHcts}). 

We continue to assume that the hypotheses of Theorem \ref{corollary:whenfinite} hold and that 
$q$ is a fixed integer with $\pi_q(X_i)$ a profinite group for each $i$. Also, let $p \geq 0$. 
The proof of Theorem \ref{corollary:whenfinite} gives that 
\[\lim_i H_p(G/N_i, \pi_q(X_i)) \cong E_2^{p,q} \cong H_p\Bigl[\lim_i (\pi_q(X_i)[(G/N_i)^\ast])\Bigr].\] Let $l \geq 0$. By 
\cite[Proposition 5.5.3]{Ribes}, there are 
the following isomorphisms of abelian groups:
\begin{align*}
\lim_i \pi_q(X_i)[(G/N_i)^l] & \cong \lim_i \oplus_{(G/N_i)^l} (\widehat{\mathbb{Z}}[[G]] 
\widehat{\otimes}_{\widehat{\mathbb{Z}}[[G]]} \pi_q(X_i)) 
\\ & 
\cong \lim_i \bigl((\oplus_{(G/N_i)^l} \widehat{\mathbb{Z}}[[G]])
\widehat{\otimes}_{\widehat{\mathbb{Z}}[[G]]} \pi_q(X_i)\bigr).\end{align*}
In the last expression, $\oplus_{(G/N_i)^l} \widehat{\mathbb{Z}}[[G]]$ is a 
profinite right $\widehat{\mathbb{Z}}[[G]]$-module and, under the $\widehat{\mathbb{Z}}[[G]]$-action on this right $\widehat{\mathbb{Z}}[[G]]$-module, $G$ acts trivially on the indexing set $(G/N_i)^l$. Then we have the isomorphisms
\begin{align*}
\lim_i \bigl((\oplus_{(G/N_i)^l} \widehat{\mathbb{Z}}[[G]])
\widehat{\otimes}_{\widehat{\mathbb{Z}}[[G]]} \pi_q(X_i)\bigr)
& \cong \lim_{j \in \{i\}} \lim_{j' \in \{i\}} 
\bigl((\oplus_{(G/N_j)^l} \widehat{\mathbb{Z}}[[G]])
\widehat{\otimes}_{\widehat{\mathbb{Z}}[[G]]} \pi_q(X_{j'})\bigr)\\
& \cong \bigl(\lim_{j \in \{i\}} (\oplus_{(G/N_j)^l} \widehat{\mathbb{Z}}[[G]])\bigr)
\widehat{\otimes}_{\widehat{\mathbb{Z}}[[G]]}\pi_q(X)\end{align*} of abelian groups, where the 
first isomorphism is by cofinality and the second one is by \cite[Lemma 5.5.2]{Ribes}. Also, 
$\lim_{j \in \{i\}} (G/N_j)^l \cong G^l$, so that by \cite[Proposition 5.2.2 and proof of 
Proposition 5.5.3, (e)]{Ribes}, $\lim_{j \in \{i\}} (\oplus_{(G/N_j)^l} \widehat{\mathbb{Z}}[[G]])$ is a free 
profinite right $\widehat{\mathbb{Z}}[[G]]$-module on the profinite space $G^l$. Here, when 
$l=0$, $G^l = \{e\}$ and $\lim_{j \in \{i\}} (\oplus_{(G/N_j)^l} \widehat{\mathbb{Z}}[[G]])$ is 
$\widehat{\mathbb{Z}}[[G]]$. 
To summarize, we see that $E_2^{p,q}$ 
is the $p$th homology of a chain complex that in degree $l$, for $l \geq 0$, is an abelian group that is isomorphic to 
\[\bigl(\lim_{j \in \{i\}} (\oplus_{(G/N_j)^l} \widehat{\mathbb{Z}}[[G]])\bigr)
\widehat{\otimes}_{\widehat{\mathbb{Z}}[[G]]}\pi_q(X),\] and the left term in this complete tensor product is a free profinite right $\widehat{\mathbb{Z}}[[G]]$-module on $G^l$. 

For the continuous homology of $G$ with coefficients in profinite right $\widehat{\mathbb{Z}}[[G]]$-modules, 
as considered in \cite[Chapter 6]{Ribes} -- whereas we use coefficients that are in profinite left $\widehat{\mathbb{Z}}[[G]]$-modules, 
one regards $\widehat{\mathbb{Z}}$ as a profinite left $\widehat{\mathbb{Z}}[[G]]$-module with 
trivial left $G$-action. As explained in \cite[page 206]{Ribes}, there is 
the inhomogeneous bar resolution of $\widehat{\mathbb{Z}}$, which is a free resolution of $\widehat{\mathbb{Z}}$ that in each degree $l$ is ${\tilde{L}}_l$, the 
free profinite left $\widehat{\mathbb{Z}}[[G]]$-module on the profinite space $G^l$. If $A$ is a profinite right $\widehat{\mathbb{Z}}[[G]]$-module, then its continuous group homology is the 
homology of a complex that in each degree $l$ is $A \widehat{\otimes}_{\widehat{\mathbb{Z}}[[G]]} {\tilde{L}}_l$. Thus, there are parallels between the chain complex whose homology gives $E_2^{p,q}$ and the computation of continuous group homology with the inhomogeneous bar resolution of $\widehat{\mathbb{Z}}$. These parallels suggest that $E_2^{p,q}$ is $H^c_p(G,\pi_q(X))$ and the 
isomorphism in (\ref{isowithHcts}) confirms that this is the case. 

Now we give the result that follows from Theorem \ref{corollary:whenfinite} by applying (\ref{isowithHcts}). 

\begin{Cor}\label{profinitecase}
Let $G$ be a profinite group and let $X = \holim_i X_i$ be an $S[[G]]$-module. Suppose that 
for every integer $t$, the inverse system $\{\pi_t(X_i)\}_i$ consists of profinite groups and continuous 
homomorphisms. Also, suppose that for all $t$ and each $i$, the induced action of the discrete group $G/N_i$ on $\pi_t(X_i)$ 
is continuous. Then the homotopy orbit spectral sequence has the form
\[E_2^{p,q} \cong H^c_p(G, \pi_q(X)) \Longrightarrow \pi_{p+q}(X_{hG}).\]
\end{Cor}

The following result -- an easy consequence of the proofs of Theorems \ref{compactHaus} and \ref{corollary:whenfinite} -- is a distillation of the key steps in these two proofs.

\begin{Thm}\label{2bigconditions}
Let $G$ be a profinite group and suppose that $X = \holim_i X_i$ is an $S[[G]]$-module. If 
\[\lims_i \pi_q(X_i)[(G/N_i)^l] = 0, \ \ \ \text{for all} \ s > 0, \ q \in \mathbb{Z}, \ \text{and} \ l \geq 0,\] 
and there is an isomorphism
\[H_p\Bigl[\lim_i (\pi_q(X_i)[(G/N_i)^\ast])\Bigr] 
\cong \lim_i H_p(\pi_q(X_i)[(G/N_i)^\ast]), \ \ \ \text{for all} \ p \geq 0, \ q \in \mathbb{Z},\] then there is a homotopy 
orbit spectral sequence of the form
\[E_2^{p,q} \cong \lim_i H_p(G/N_i, \pi_q(X_i)) \Longrightarrow \pi_{p+q}(X_{hG}).\] 
\end{Thm}

\section{The homotopy orbit spectral sequence 
when $G$ is countably based}\label{tools}

\par
Now we study the homotopy orbit spectral 
sequence of Theorem \ref{generalh-oSS} in the case when 
the profinite group $G$ is countably based. In the following remark, we lay out a homological construction, a special case of which 
is used in our next result, Theorem \ref{oss}.

\begin{Rk}\label{limonegivescomplex}
Let $J$ be a directed poset, and let $\{A_\ast^j\}$ be a $J$-indexed inverse system of 
chain complexes $A_\ast^j$ in $\mathrm{Ch}_{\geq 0}$, the 
category of non-negatively graded chain complexes in $\mathbf{Ab}$ (the category of abelian groups): 
$\{A_\ast^j\}$ is an object in 
$(\mathrm{Ch}_{\geq 0})^{J^\mathrm{op}}$, the category of functors 
$J^\mathrm{op} \to \mathrm{Ch}_{\geq 0}$. For each $j \in J$, $d^j_\ast$ denotes the differentials of $A_\ast^j$: for every $k \geq 1$, each $d^j_k$ is a homomorphism $A^j_k \to 
A^j_{k-1}$.  Notice that for $s \geq 0$, 
the fact that the functor 
$\displaystyle{\lims_J (-):\textbf{Ab}^{J^\mathrm{op}} \to \textbf{Ab}}$ is 
additive means that whenever 
$k \geq 1$, for the 
morphisms $\{d^j_{k+1}\}$, $\{d_k^j\}$ in $\mathbf{Ab}^{J^\mathrm{op}}$, we have
\begin{align*} 
(\lims_J \{d^j_k\}) \circ (\lims_J \{d^j_{k+1}\}) & = 
\lims_J \{d_k^j\circ d_{k+1}^j \} = \lims_J \{0 \: A^j_{k+1} \to A^j_{k-1} \} \\
& = (0 \: \lims_j A^j_{k+1} \to \lims_j A^j_{k-1}).   
\end{align*}  
It follows that for every $s \geq 0$, 
$\displaystyle{\lims_J (-)}$ extends to a 
functor \[\lims_J (-) \: 
(\mathrm{Ch}_{\geq 0})^{J^\mathrm{op}} \to \mathrm{Ch}_{\geq 0}, \ \ \ 
\{A_\ast^j\} \ \mapsto \ \lims_J \{A_\ast^j\} = \lims_j A_\ast^j,
\] where for each $k \geq 0$, 
\[(\lims_j A_\ast^j)_k := \lims_j A^j_k\] defines the 
group of chains of the chain complex 
$\lim^s_j A_\ast^j$ in degree $k$. Thus, for each $p, s \geq 0$, one 
can form \[H_p\mspace{-5mu}\left[\lims_j A^j_\ast\right] = H_p\mspace{-5mu}\left[\cdots \to 
\lims_j A^j_k \to \cdots \to \lims_j A^j_1 \to \lims_j A^j_0\right],\]  
the $p$th homology of the complex $\lim^s_j A_\ast^j$.
\end{Rk}

To help with understanding 
our next result, we recall the following. If $G$ is any profinite group and $\holim_i X_i$ is an $S[[G]]$-module, then for 
each $i$, there is the simplicial spectrum 
$(X_i[(G/N_i)^\bullet])_\mathtt{f}$, so that as in the proof of Theorem \ref{corollary:whenfinite}, for every integer $q$, there is an 
isomorphism 
\begin{equation}\label{iso-of-complexes}\zig
\bigl\{\pi_q((X_i[(G/N_i)^\ast])_\mathtt{f})\bigr\} 
\cong \bigl\{\pi_q(X_i)[(G/N_i)^\ast]\bigr\}
\end{equation}
in $(\mathrm{Ch}_{\geq 0})^{(\{i\}^\mathrm{op})}\mspace{-2mu}$, 
the category of $\{i\}$-indexed inverse systems of non-negatively graded chain complexes. In 
(\ref{iso-of-complexes}), the left-hand side is the diagram of 
Moore complexes and on the right-hand side, each complex 
$\pi_q(X_i)[(G/N_i)^\ast]$ is defined as in Definition \ref{pixgbullet}. 

\begin{Rk}
We define some helpful notation. When an exact sequence in $\mathbf{Ab}$ of the form 
$A \to B \to C$ extends beyond a single line, we write it as 
\begin{align*} A & \to B \scriptstyle{\ - \cdots}\\ 
& \  
\to C.\end{align*}\end{Rk}

We recall from Definition \ref{typeofmodule} and Remark \ref{rkforbased} that ``let $\holim_i X_i$ be 
a countably based $S[[G]]$-module" means that $G$ is a countably based profinite group, 
there is a fixed descending chain $\{N_i\}_{i \geq 0}$ of open normal subgroups of $G$ with 
$G \cong \lim_{i \geq 0} G/N_i$, and there is a tower $\{X_i\}_{i \geq 0}$ of $G$-spectra and 
$G$-equivariant maps, such that the pair $(\{X_i\}_{i \geq 0}, \holim_i X_i)$ is an $S[[G]]$-module.

\begin{Thm}\label{oss}
Let $X = \holim_i X_i$ be a countably based $S[[G]]$-module. 
For the $E_2$-term of the homotopy orbit spectral sequence 
\[E_2^{p,q} = H_p(\pi_q(\holim_i (X_i[(G/N_i)^\ast])_\mathtt{f})) 
\Longrightarrow \pi_{p+q}(X_{hG})\] and the tower 
$\bigl\{\pi_q(X_i)[(G/N_i)^\ast]\bigr\}$ of chain complexes, where $q$ is any integer, 
there is a long exact sequence 
\begin{align*}  
& \cdots \xrightarrow{\,\scriptscriptstyle{\partial}\,} 
H_p\mspace{-5mu}\left[\limone_i \pi_{q+1}(X_i)[(G/N_i)^\ast]\right] \to E_2^{p,q}\to H_p\mspace{-5mu}\left[\lim_i \pi_{q}(X_i)[(G/N_i)^\ast]\right] 
\scriptstyle{\ - \cdots}\\
& \ \ \ \xrightarrow{\,\scriptscriptstyle{\partial}\,} 
H_{p-1}\mspace{-5mu}\left[\limone_i \pi_{q+1}(X_i)[(G/N_i)^\ast]\right] \to \cdots \to H_1\mspace{-5mu}\left[\lim_i \pi_{q}(X_i)[(G/N_i)^\ast]\right] 
\scriptstyle{\ - \cdots}\\
& \ \ \ \ \ \xrightarrow{\,\scriptscriptstyle{\partial}\,} 
H_0\mspace{-5mu}\left[\limone_i \pi_{q+1}(X_i)[(G/N_i)^\ast]\right] \to E_2^{0,q}\to H_0\mspace{-5mu}\left[\lim_i \pi_{q}(X_i)[(G/N_i)^\ast]\right] \xrightarrow{\,\scriptscriptstyle{\partial}\,} 0. 
\end{align*}
\end{Thm}

\begin{proof}
Let $q$ be any integer. For each $l \geq 0$, there is the 
Milnor short exact sequence
\[
0 \mspace{-2mu}\shortrightarrow\mspace{-2mu} \limone_i \pi_{q+1}(X_i)[(G/N_i)^l] \mspace{-2mu}\shortrightarrow\mspace{-2mu}
\pi_q(\holim_i (X_i[(G/N_i)^l])_\mathtt{f}) 
\mspace{-2mu}\shortrightarrow\mspace{-2mu} \lim_i \pi_q(X_i)[(G/N_i)^l] \mspace{-2mu}\shortrightarrow\mspace{-2mu} 0.
\] 
As explained in Remark \ref{limonegivescomplex}, $\displaystyle{\limone_i \pi_{q+1}(X_i)[(G/N_i)^\ast]}$ is a chain complex, so that by letting $l$ vary, the above 
Milnor short exact sequences 
give the short exact sequence 
\[
0 \mspace{-2mu}\shortrightarrow\mspace{-2mu} \limone_i \pi_{q+1}(X_i)[(G/N_i)^\ast] \mspace{-2mu}\shortrightarrow\mspace{-2mu}
\pi_q(\holim_i (X_i[(G/N_i)^\ast])_\mathtt{f}) 
\mspace{-2mu}\shortrightarrow\mspace{-2mu} \lim_i \pi_q(X_i)[(G/N_i)^\ast] \mspace{-2mu}\shortrightarrow\mspace{-2mu} 0
\] 
of chain complexes, and associated to this last short exact sequence 
is the desired long exact sequence of 
homology groups. 
\end{proof}

Now we show that the ``term on the left" in ``the degree zero row" of 
the long exact sequence of Theorem \ref{oss} can be simplified. We 
use the standard notation that if $K$ is an abstract group, then
\[(-)_{\mspace{-1.5mu}{\scriptscriptstyle{K}}} \: _{\mathbb{Z}[K]}\mathrm{Mod} \to \mathbf{Ab}, \ \ \ M 
\mapsto M_{\mspace{-1.5mu}{\scriptscriptstyle{K}}}\]
is the right exact coinvariants functor from $K$-modules to abelian groups whose left derived functors are group homology. Thus, for any $K$-module $M$, there is an isomorphism 
$H_0(K, M) \cong M_{\mspace{-1.5mu}{\scriptscriptstyle{K}}}$. 

\begin{Thm}\label{Hzero}
Suppose that $\holim_i X_i$ is a countably based $S[[G]]$-module. 
Then there is an isomorphism 
\[H_0\mspace{-5mu}\left[\limone_i \pi_{q+1}(X_i)[(G/N_i)^\ast]\right] 
\cong \limone_i (\pi_{q+1}(X_i))_{\mspace{-1mu}{\scriptscriptstyle{G/N_i}}}\mspace{1mu},\] where $q$ is 
any integer and the right-hand side is $\lim^1_i(-)$ applied to a tower of various coinvariants.
\end{Thm}

\begin{proof}
For each $i \geq 0$, there is the commutative diagram 
\[\xymatrix{
\pi_{q+1}(X_{i+1})[G/N_{i+1}] \ar[d] \ar^-{d^{i+1}_1}[r] & \pi_{q+1}(X_{i+1}) 
\ar[d] \ar^-{\pi^{i+1}}[r] & \pi_{q+1}(X_{i+1})/{\mathrm{im}(d^{i+1}_1)} \ar[d] \ar[r] & 0 \ar[d]
\\
\pi_{q+1}(X_{i})[G/N_i] \ar^-{d^{i}_1}[r] 
& \pi_{q+1}(X_i) \ar^-{\pi^{i}}[r] 
& \pi_{q+1}(X_i)/{\mathrm{im}(d^{i}_1)} \ar[r] & 0\rlap{$\mspace{2mu}$,}
}\] with both rows exact, for the following reasons: (a) the commutative 
square on the left is just a piece of the tower of chain complexes 
$\{\pi_{q+1}(X_j)[(G/N_j)^\ast]\}_{j \geq 0}$, showing the first differential for the 
complex at heights $i$ and $i+1$; (b) the homomorphisms $\pi^i$ and 
$\pi^{i+1}$ are the canonical maps to the respective cokernels; and 
(c) the map 
between these two cokernels is the unique map induced by the first two 
vertical maps on the left (this uses the exactness of the rows). Letting 
$i$ vary implies that there is the diagram 
\begin{equation}\label{indexedby4}\zig
\{\pi_{q+1}(X_{i})[G/N_i]\}
\xrightarrow{\,\{d^i_1\}\,} \{\pi_{q+1}(X_i)\}
\xrightarrow{\,\{\pi^i\}\,} \{\pi_{q+1}(X_i)/{\mathrm{im}(d^{i}_1)}\}
\to \{0\}
,
\end{equation} 
with exact rows, 
in the category $\mathbf{tow}(\mathbf{Ab})$ of 
towers of abelian 
groups.

Now suppose that 
\[\{0\} \to \{A_i\} \to \{B_i\} \to \{C_i\} \to \{0\}\] is a short exact sequence 
in $\mathbf{tow}(\mathbf{Ab})$. Then the sequence 
\[0 \to \lim_i A_i \to \lim_i B_i \to \lim_i C_i \to \limone_i A_i \to 
\limone_i B_i \to \limone_i C_i \to 0\] is exact, so that the additive functor 
$\lim^1_i \: \mathbf{tow}(\mathbf{Ab}) \to \mathbf{Ab}$ is right exact. 
Thus, applying $\lim^1_i (-)$ to diagram (\ref{indexedby4}) gives 
the exact sequence
\[\limone_i \pi_{q+1}(X_{i})[G/N_i] \xrightarrow{\,\lim^1_i d^i_1\,} 
\limone_i \pi_{q+1}(X_i)  
\xrightarrow{\,\lim^1_i \pi^i\,} \limone_i \pi_{q+1}(X_i)/{\mathrm{im}(d^{i}_1)} \to 0,\] and hence, 
\begin{align*}
H_0\mspace{-5mu}\left[\limone_i \pi_{q+1}(X_i)[(G/N_i)^\ast]\right] 
& = \bigl(\limone_i \pi_{q+1}(X_i)\bigr)/{\mathrm{im}(\limone_i d^i_1)} 
\\ & \cong \limone_i \pi_{q+1}(X_i)/{\mathrm{im}(d^{i}_1)} \\ 
& = \limone_i H_0(\pi_{q+1}(X_i)[(G/N_i)^\ast]) \\ &\cong \limone_i H_0(G/N_i, \pi_{q+1}(X_i)),
\end{align*} as desired.
\end{proof}

Our next result can help with computing the homology groups 
\[H_\ast\mspace{-5mu}\left[\limone_i \pi_{q+1}(X_i)[(G/N_i)^\ast]\right]\] in Theorem 
\ref{oss}. 

\begin{Thm}\label{limlimone}
For any $q \in \mathbb{Z}$ and each $l \geq 0$, the $l$th group of chains in the chain complex 
$\lim^1_i \pi_{q}(X_i)[(G/N_i)^\ast]$ that appears in $\mathrm{Theorem}$ $\mathrm{\ref{oss}}$ satisfies the 
isomorphism 
\[\limone_i \pi_{q}(X_i)[(G/N_i)^l] \cong  \lim_j \bigl((\limone_i \pi_{q}(X_i))[(G/N_j)^l]\bigr),\] 
where the limit $\,\displaystyle{\lim_j}$ above is indexed by $\{j\} = \{i\}$ and, for each $j$, the expression 
$(\displaystyle{\limone_i \pi_{q}(X_i))[(G/N_j)^l]}$ on the right-hand side is 
$\textstyle{\bigoplus_{(G/N_j)^l}} (\displaystyle{\limone_i \pi_{q}(X_i))}$. 
\end{Thm}

\begin{proof}
Let $l \geq 0$ be fixed. By \cite[Lemma 1.10]{AndersonPhantom}, there is the identity
\[\limone_i \pi_{q}(X_i)[(G/N_i)^l] = \limone_{j,i} \pi_{q}(X_i)[(G/N_j)^l],\] where the last 
expression is the first derived functor of limit for double towers. Then by \cite[page 429]{Caruso}, 
there is a short exact sequence
\begin{align*}
0 \to \limone_j \bigl(\lim_i (\pi_{q}(X_i)[(G/N_j)^l])\bigr) & \to \limone_i \pi_{q}(X_i)[(G/N_i)^l] 
\scriptstyle{\ - \cdots}\\ & \ \ \ \ \ \ \to \lim_j \bigl(\limone_i (\pi_{q}(X_i)[(G/N_j)^l])\bigr) \to 0.
\end{align*} Notice that there are isomorphisms
\begin{align*}
\limone_j & \bigl(\lim_i (\pi_{q}(X_i)[(G/N_j)^l])\bigr) \cong \limone_j \bigl(\lim_i \textstyle{\prod}_{(G/N_j)^l} \pi_{q}(X_i)\bigr)\\ & \cong \limone_j \textstyle{\prod}_{(G/N_j)^l} \displaystyle{\lim_i \pi_{q}(X_i)} \cong \limone_j (\lim_i \pi_{q}(X_i))[(G/N_j)^l] = 0,\end{align*} 
where the last step 
applies the fact that the penultimate expression is $\lim^1_j$ applied to a tower of surjections. 
Then the short exact sequence yields 
\begin{align*}
\limone_i & \pi_{q}(X_i)[(G/N_i)^l] \cong
\lim_j \bigl(\limone_i (\pi_{q}(X_i)[(G/N_j)^l])\bigr)\\ & \cong \lim_j \textstyle{\coprod}_{(G/N_j)^l} \displaystyle{\limone_i \pi_{q}(X_i)} \cong \lim_j \bigl((\limone_i \pi_{q}(X_i))[(G/N_j)^l]\bigr),
\end{align*} where the second isomorphism follows from the fact that 
the functor $\lim^1_i(-)$ is additive and for each fixed $j$, 
the tower $\{\pi_{q}(X_i)[(G/N_j)^l]\}_i$ is the finite coproduct $\textstyle{\coprod}_{(G/N_j)^l} 
\{\pi_{q}(X_i)\}_i$ in the functor category $\mathbf{Ab}^{(\{i\}^\mathrm{op})}$.
\end{proof}

\begin{Rk}\label{quickRk}
Given an arbitrary countably based $S[[G]]$-module $\holim_i X_i$, it follows immediately from Theorem \ref{oss} that if there is some integer $q$ for which 
$\displaystyle{\limone_i  \pi_{q+1}(X_i)[(G/N_i)^l]} = 0,$ for all $l \geq 0$, then 
for 
every $p \geq 0$, there is an isomorphism
$E_2^{p,q} \xrightarrow{\,\cong\,} H_p\mspace{-1mu}\Bigl[\displaystyle{\lim_i  \pi_{q}(X_i)[(G/N_i)^\ast]}\Bigr].$
\end{Rk}
The following result is a straightforward consequence of Theorem \ref{limlimone} and Remark \ref{quickRk}. 

\begin{Cor}\label{mightneed}
Let $\holim_i X_i$ be a countably based $S[[G]]$-module. If $q$ is an integer such that 
$\displaystyle{\limone_i \pi_{q+1}(X_i) = 0}$, then for the $E_2$-term of the homotopy orbit spectral sequence, there is an isomorphism
\[E_2^{p,q} \xrightarrow{\,\cong\,} H_p\mspace{-1mu}\Bigl[\displaystyle{\lim_i  \pi_{q}(X_i)[(G/N_i)^\ast]}\Bigr],\] for each $p \geq 0$.
\end{Cor}

\begin{Thm}\label{specialcase}
Let $q$ be any integer and suppose that $X = \holim_i X_i$ is a countably based $S[[G]]$-module. 
If the towers $\{\pi_{q}(X_i)\}_i$ and $\{\pi_{q+1}(X_i)\}_i$ are diagrams in the category of compact Hausdorff abelian groups, and for every $i$, the action of the discrete group $G/N_i$ on $\pi_q(X_i)$ is continuous, then the $E_2$-term of the 
homotopy orbit spectral sequence satisfies the isomorphism 
\[E_2^{\ast,q} \cong \lim_i H_\ast(G/N_i, \pi_q(X_i))\] of non-negatively graded abelian groups. 
\end{Thm}
\begin{proof}
Notice that $\lim^1_i \pi_{q+1}(X_i) = 0$. Then Corollary \ref{mightneed} yields an isomorphism 
$E_2^{\ast,q} \cong H_\ast\mspace{-3mu}\Bigl[\lim_i  \pi_{q}(X_i)[(G/N_i)^\ast]\Bigr]$ of non-negatively graded abelian groups. As in the proof of Theorem \ref{corollary:whenfinite}, the tower $\{\pi_q(X_i)[(G/N_i)^\ast]\}_i$ of chain complexes lives in $\mathcal{CHA}$. The argument in the last paragraph of the proof of Theorem \ref{corollary:whenfinite} completes the proof. 
\end{proof}
%
%

The next result is a consequence of Corollary \ref{mightneed}, \cite[Theorem 3.5.8]{Weibel}, and (\ref{homology}).

\begin{Thm}\label{gottenwithML} 
Let $q$ be any integer and let $X = \holim_i X_i$ be a countably based $S[[G]]$-module. Also, 
let $E_2^{\ast,\ast}$ denote the $E_2$-term of its homotopy orbit spectral sequence. 
If $\lim^1_i \pi_{q+1}(X_i) = 0$ and, for each $l \geq 0$, the tower 
$\{\pi_{q}(X_i)[(G/N_i)^l]\}_i$ 
satisfies the Mittag-Leffler condition, 
then there is the short exact sequence 
\[0 \rightarrow \limone_i H_{\ast+1}(G/N_i,\pi_q(X_i)) \rightarrow 
E_2^{\ast,q} \rightarrow \lim_i H_\ast(G/N_i,\pi_q(X_i)) \rightarrow 0\] 
of non-negatively graded abelian groups. 
\end{Thm}

\begin{Rk}\label{finitelygenerated}
By \cite[Corollary 6.5.10]{Weibel}, when $K$ is a finite group and $M$ is a finitely generated $K$-module, $H_p(K,M)$ is finite for all $p > 0$. Thus, when the hypotheses of Theorem \ref{gottenwithML} hold and additionally, for each $i$, 
$\pi_q(X_i)$ is a finitely generated $G/N_i$-module, the map
$E_2^{\ast,q} \rightarrow \lim_i H_\ast(G/N_i,\pi_q(X_i))$ is an isomorphism.
\end{Rk} 

\end{section}
\begin{section}{Eilenberg-Mac Lane spectra and their homotopy orbits}\label{EM}
\par

Let $G$ be any profinite group and, as in Definition \ref{nicesystem}, let 
$\{A_i\}$ be a nice inverse system of $G$-modules with respect to the 
collection $\{N_i\}$ of open normal subgroups. 

Let $\Gamma \: \mathbf{Ch}_+ \rightarrow \mathtt{s}(\mathbf{Ab})$ be the functor in the 
Dold-Kan correspondence from $\mathbf{Ch}_+$, the category of chain complexes 
$C_\ast$ with $C_n=0$ for $n <0$, to $\mathtt{s}(\mathbf{Ab})$, the category of simplicial 
abelian groups (see, for example, \cite[Chapter III, Corollary 2.3]{GJ}). 
Also, if $A$ is an abelian group, 
let $A[-n]$ be the chain complex that is $A$ in degree $n$ and 
zero elsewhere.
\par
Given the inverse system $\{A_i\}$, we explain how to form the inverse system 
$\{H(A_i)\}$ of Eilenberg-Mac Lane spectra, by following the construction given in \cite{Jardinesummary}. Then the pair $(\{H(A_i)\}, 
\holim_i H(A_i))$ is an 
$S[[G]]$-module. By functoriality, for each $k \geq 0$, 
$\{\Gamma(A_i[-k])\}$ is an inverse system of 
simplicial $G$-modules and $G$-equivariant maps, such that, for each $i$, 
$\Gamma(A_i[-k])$ is the Eilenberg-Mac Lane space $K(A_i,k)$ and 
$\Gamma(A_i[-k])$ is a simplicial $G/N_i$-module. Furthermore, by taking $0$ as 
the basepoint, each $\Gamma(A_i[-k])$ is a pointed simplicial set. 
\par
For each $i$, we define the Eilenberg-Mac Lane spectrum $H(A_i)$ by 
$(H(A_i))_k = \Gamma(A_i[-k]),$ so that $\pi_0(H(A_i)) = A_i$ and $\pi_n(H(A_i)) = 0$, 
when $n \neq 0$. Then, by functoriality, 
$\{H(A_i)\}$ is an inverse system of 
$G$-spectra and $G$-equivariant maps, such that each $H(A_i)$ is a 
$G/N_i$-spectrum. Since each $(H(A_i))_k$ is a fibrant simplicial set and each 
$H(A_i)$ is an 
$\Omega$-spectrum (see, for example, \cite[Example 21]{Jardinesummary}), $H(A_i)$ 
is a fibrant spectrum. These facts imply the following result.
\begin{Lem}
If $G$ is any profinite group and $\{A_i\}$ is a nice inverse 
system of $G$-modules with respect to $\{N_i\}$, then the pair $(\{H(A_i)\}, \holim_i H(A_i))$ is an $S[[G]]$-module. 
\end{Lem} 

By Theorem \ref{generalh-oSS}, there is the homotopy orbit spectral sequence
\[E_2^{p,q} \Longrightarrow \pi_{p+q}((\holim_i H(A_i))_{hG}).\] For each $l \geq 0$, 
there is the homotopy spectral sequence 
\[^l\mspace{-2mu}E^{s,t}_2  \Longrightarrow \pi_{t-s}(\holim_i (H(A_i)[(G/N_i)^l])_\mathtt{f}),\] where
\begin{align*}
^l\mspace{-2mu}E_2^{s,t} & = \lims_i \pi_t((H(A_i)[(G/N_i)^l])_\mathtt{f}) \cong \lims_i \pi_t(H(A_i))[(G/N_i)^l]\\
& = 
\begin{cases}
0, & \text{if} \ t \neq 0;\\ 
\displaystyle{\lim_i}^s A_i[(G/N_i)^l], & \text{if} \ t = 0,\end{cases}\end{align*} so that this homotopy spectral sequence collapses, giving 
\[\pi_q(\holim_i (H(A_i)[(G/N_i)^l])_\mathtt{f}) \cong \displaystyle{\lim_i}^q A_i[(G/N_i)^l], \ \ \ q \in \mathbb{Z}.\] In this last isomorphism, for each $q < 0$, the right-hand side is $0$. Using 
the notation of Remark \ref{limonegivescomplex}, we conclude that when $G$ is any profinite group and $\{A_i\}$ is a nice inverse system of $G$-modules for $\{N_i\}$, then the 
$E_2$-term of the homotopy orbit spectral sequence has the form 
\begin{equation}\label{e2term}\zig
E_2^{p,q} \cong 
\begin{cases}
H_p \mspace{-3.5mu}\left[\displaystyle{\lim_i}^q A_i[(G/N_i)^\ast]\right], & q \geq 0;\\
0, & q < 0,\end{cases}\end{equation} where 
\[A_i[(G/N_i)^\ast] := \pi_0(H(A_i))[(G/N_i)^\ast], \ \ \ \text{for each} \ i,\] is 
defined as in Definition 
\ref{pixgbullet}.  

Below, if $K$ is a finite group and $M$ is a $K$-module, then $H_p(K,M) = 0$, 
whenever $p < 0$. 

\begin{Thm}\label{generalE-Mresult}
Let $G$ be a profinite group and let $\{A_i\}$ be a nice inverse 
system of $G$-modules with respect to $\{N_i\}$. If the $E_2$-term of the 
homotopy orbit spectral sequence for 
the $S[[G]]$-module $\holim_i H(A_i)$ satisfies the isomorphism 
$E_2^{p,q} \cong \lim_i H_p(G/N_i, \pi_q(H(A_i)))$ for 
all $p \geq 0$, $q \in \mathbb{Z}$, then 
there is the $\mathbb{Z}$-graded isomorphism
\[\pi_\ast((\holim_i H(A_i))_{hG}) \cong \lim_i H_\ast(G/N_i, A_i).\] 
\end{Thm} 

\begin{proof}
The isomorphism satisfied by the $E_2$-term implies that
\[E_2^{p,q} \cong \begin{cases} 0, & 
\mathrm{if} \ q \neq 0; \\ \displaystyle{\lim_i H_p(G/N_i, A_i)}, & \mathrm{if} \ q = 0,\end{cases}\] 
so that the homotopy orbit spectral sequence collapses, giving the desired conclusion.
\end{proof}

It is straightforward to see that Theorem \ref{corollary:whenfinite}, Corollary \ref{profinitecase}, Theorem \ref{gottenwithML}, and Remark \ref{finitelygenerated} give conditions that result in the hypotheses of Theorem \ref{generalE-Mresult} holding. For example, there 
is the following result, whose statement makes use of the inverse system 
$\{A_i[(G/N_i)^l]\}$, for $l \geq 0$, that plays a role in (\ref{e2term}). 

\begin{Cor}\label{E-Mf-g}
Let $G$ be a countably based profinite group, with $\{N_i\}$ a descending chain of open normal subgroups in the 
sense of $\mathrm{(\ref{chainofopens})}$, and let $\{A_i\}$ be a nice inverse 
system of $G$-modules with respect to $\{N_i\}$. If the tower $\{A_i[(G/N_i)^l]\}$ of abelian groups satisfies the Mittag-Leffler condition for each $l \geq 0$ and, for each $i$, $A_i$ is a finitely generated $G/N_i$-module, then there is the $\mathbb{Z}$-graded isomorphism
\[\pi_\ast((\holim_i H(A_i))_{hG}) \cong \lim_i H_\ast(G/N_i, A_i).\] 
\end{Cor}

\begin{proof}
Let $q$ be an integer. The tower $\{\pi_{q+1}(H(A_i))\}$ is either 
$\{\pi_0(H(A_i))\} \cong \{A_i\}$, which satisfies the Mittag-Leffler condition, 
or the tower $\{0\}$, which is the trivial group in each level, and so for all $q$, $\lim^1_i \pi_{q+1}(H(A_i)) = 0$. Also, given any integer $q$ and any $l \geq 0$, the tower $\{\pi_q(H(A_i))[(G/N_i)^l]\}$ is 
either $\{\pi_0(H(A_i))[(G/N_i)^l]\} \cong \{A_i[(G/N_i)^l]\}$ or $\{0\}$, both of 
which satisfy the Mittag-Leffler condition. 

Now let $p > 0$ and let $i$ be arbitrary: as in Remark \ref{finitelygenerated}, 
\[H_p(G/N_i, \pi_0(H(A_i))) \cong H_p(G/N_i, A_i)\] is finite, and if $q \neq 0$, 
$H_p(G/N_i, \pi_q(H(A_i))) = 0$. Thus, for every $p > 0$ and 
each $q$, $\lim^1_i H_p(G/N_i, \pi_q(H(A_i))) = 0$. These observations, 
together with Theorem \ref{gottenwithML}, imply that Theorem \ref{generalE-Mresult} applies.
\end{proof}

\end{section}

\end{document}